\documentclass[twocolumn]{autart} 

\usepackage{epsfig} 
\usepackage{amsmath} 
\usepackage{amssymb} 
\usepackage{euscript}
\usepackage{graphicx}       
\usepackage{color}       
\newcommand{\bfsfb}{\mathbf{\mathsf{b}}}
\newcommand{\bfsfc}{\mathbf{\mathsf{c}}}

\newfont{\Bb}{msbm10 scaled\magstep0}
\def\IR{\mbox {\Bb R}}
\def\IC{\mbox {\Bb C}}

\newcommand{\bfA}{ {\mathbf A} }
\newcommand{\bfB}{ {\mathbf B} }
\newcommand{\bfC}{ {\mathbf C} }

\newcommand{\bfE}{ {\mathbf E} }
\newcommand{\bfF}{ {\mathbf F} }
\newcommand{\bfG}{ {\mathbf G} }
\newcommand{\bfH}{ {\mathbf H} }
\newcommand{\bfI}{ {\mathbf I} }
\newcommand{\bfJ}{ {\mathbf J} }
\newcommand{\bfK}{ {\mathbf K} }
\newcommand{\bfL}{ {\mathbf L} }
\newcommand{\bfM}{ {\mathbf M} }

\newcommand{\bfQ}{ {\mathbf Q} }
\newcommand{\bfR}{ {\mathbf R} }
\newcommand{\bfS}{ {\mathbf S} }
\newcommand{\bfT}{ {\mathbf T} }

\newcommand{\bfV}{ {\mathbf V} }
\newcommand{\bfW}{ {\mathbf W} }
\newcommand{\bfX}{ {\mathbf X} }

\newcommand{\bfVr}{{\bfV}_r }

\newcommand{\hbfVr}{ \widehat{\bfV}_r }
\newcommand{\hbfWr}{ \widehat{\bfW}_r }

\newcommand{\htwo}{{ {\mathcal H}_2} }
\newcommand{\hinf}{{{\mathcal H}_{\infty}}}
\newcommand{\Hinf}{{{\mathcal H}_{\infty}}}
\newcommand{\Htwo}{{ {\mathcal H}_2}}  
\newcommand{\PH}{{ \mathcal{P\!H} }}

\newcommand{\tbfJ}{ \widetilde{\bfJ}}
\newcommand{\tbfR}{ \widetilde{\bfR}}
\newcommand{\tbfQ}{ \widetilde{\bfQ}}
\newcommand{\tbfB}{ \widetilde{\bfB}}
\newcommand{\tbfx}{ \widetilde{\bfx}}

\newcommand{\bfe}{ {\mathbf e} }

\newcommand{\bfu}{ {\mathbf u} }
\newcommand{\bfv}{ {\mathbf v} }

\newcommand{\bfx}{ {\mathbf x} }
\newcommand{\bfy}{ {\mathbf y} }
\newcommand{\bfz}{ {\mathbf z} }

\newcommand{\phss}{port-Hamiltonian systems }

\newcommand{\ph}{port-Hamiltonian }

\newtheorem{example}{Example}

\begin{document}

\begin{frontmatter}
\title{Structure-preserving tangential interpolation for model reduction of port-Hamiltonian Systems}

\thanks[footnoteinfo]{This paper has not been presented at any IFAC 
meeting. Corresponding author S.~Gugercin. Tel. +1-540-231-6549. 
Fax +1-540-231-5960}

\author[GugBea]{Serkan Gugercin}\ead{gugercin@math.vt.edu},   
\author[Polyuga]{Rostyslav V. Polyuga}\ead{\textcolor{black}{rostyslav.polyuga@gmail.com}},    
\author[GugBea]{Christopher Beattie}\ead{beattie@vt.edu},  and 
\author[vdSchaft]{Arjan van der Schaft}\ead{A.J.van.der.Schaft@rug.nl}

\address[GugBea]{Department of Mathematics, Virginia Tech, Blacksburg, VA, 24061-0123, USA}
\address[Polyuga]{ABN AMRO N.V. Bank, Gustav Mahlerlaan 10, 1082 PP (PAC HQ2015) Amsterdam, The Netherlands}
\address[vdSchaft]{Johann Bernoulli Institute for Mathematics and Computer Science, University of Groningen, P.O.Box 407,\\  9700 AK Groningen, The Netherlands}



\begin{abstract}
Port-Hamiltonian systems result from  port-based network modeling of physical systems and are an important example of passive state-space systems. In this paper, we develop
the framework for model reduction of large-scale multi-input/multi-output port-Hamiltonian systems via tangential rational interpolation. The resulting reduced-order model 
not only is a rational tangential interpolant but also retains the 
port-Hamiltonian structure; hence is passive. This reduction methodology is described in both energy and co-energy system coordinates. We also introduce an $\Htwo$-inspired algorithm for effectively choosing the interpolation points and tangential directions. The algorithm leads a reduced port-Hamiltonian model that satisfies a subset of $\Htwo$-optimality conditions.  We present several numerical examples that illustrate the effectiveness of the proposed method showing that it outperforms other existing techniques in both quality and numerical efficiency. 
\end{abstract}
\end{frontmatter}


\section{Introduction}
\label{sec:introduction}

Port-based network modeling \cite{geoplexbook2009} of physical systems exploits the common circumstance that the system under study is decomposable into (possibly many) subsystems which are interconnected through (vector) pairs of variables, whose product gives the power exchanged among subsystems. This approach is especially useful for multi-physics systems, where the subsystems may describe different categories of physical phenomena (e.g, mechanical, electrical, or hydraulic).
This approach leads
directly to \emph{port-Hamiltonian system} representations
(see \cite{Schaft00,SchaftMaschke1995,Schaft2000b,geoplexbook2009} and references therein) 
encoding structural properties related to the manner in which energy is distributed and flows through the system.   

Quite apart from underlying port-Hamiltonian structure,  models of complex physical systems often
involve discretized systems of coupled partial differential equations which lead, in turn, to dynamic models having very large state-space dimension.  
This creates a need for  \emph{model reduction} methods capable of producing 
(comparatively) low dimension surrogate models that are able to  mimic closely the original system's input/output map. 
When the system of interest additionally has port-Hamiltonian structure, then it is highly desirable for the reduced model to preserve port-Hamiltonian structure so as to maintain a variety of  system properties, such as energy
conservation, passivity, and existence of conservation laws.

 Preservation of port-Hamiltonian structure in the reduction of large scale port-Hamiltonian systems has been considered in \cite{PolvdSc10_MTNS,PolvdSc09,PolvdSc08} using Gramian-based methods; we briefly review one such approach in Section \ref{sec:effort}.   For the special case of a single-input/single-output (SISO) system, the use of (rational) Krylov methods  was employed in \cite{Polyuga_vdSchaft2009_infinity,gugercin2009ibh,Polyuga_vdSchaft2010_s0,Lohmann_Wolf_Eid_Kotyczka2009}.  In particular,  \cite{Polyuga_vdSchaft2009_infinity,Polyuga_vdSchaft2010_s0,Lohmann_Wolf_Eid_Kotyczka2009} all deal with system interpolation at single points only. 
 In this work,  we develop multi-point rational  tangential interpolation of multi-input/multi-output (MIMO) port-Hamiltonian systems.   Preservation of port-Hamiltonian structure by reducing the underlying full-order Dirac structure was presented in \cite{Polyuga_vdSchaft2011_Effort_flow,vdSchaft_Polyuga_09cdc}.
 A perturbation approach is considered in \cite{Hartmann08,Hartmann09}.   See \cite{Polyuga_2010_thesis} for an overview of recent port-Hamiltonian model reduction methods.

 For general MIMO dynamical systems, interpolatory model reduction methods produce reduced-order models whose transfer functions interpolate the original system transfer function at selected points in the complex (frequency) plane along selected input and output directions.  The main implementation cost involves solving (typically sparse) linear systems, which gives a significant advantage in large-scale settings over competing Gramian-based methods (such as balanced truncation) that must contend with a variety of large-scale dense matrix operations.  Indeed, a Schur decomposition is required for computing an \emph{exact} Gramian although recent advances using ADI-type algorithms 
have allowed the iterative computation of \emph{approximate} Gramians  for large-scale systems;  see, for example, 
\cite{penzl200clr,benner2003statespace,gugercin2003amodified,stykel2004gbm,sorensen2002thesylvester,heinkenschloss2008btm,freitas08gbr,sabino2007sls} and  references therein.

Until recently, there were no systematic strategies for selecting interpolation points and directions, but this has largely been resolved by Gugercin {\it et al.} 
\cite{gugercin2005irk,gugercin2006rki,gugercin2008hmr} where an interpolation point/tangent direction selection strategy leading to  $\Htwo$-optimal reduced order models has been proposed.
 See \cite{vandooren2008hom,bunse-gerstner2009hom,kubalinska2007h0i} for related work and \cite{antoulas2010imr} for a recent survey.

The goal of this work is to demonstrate that interpolatory model reduction techniques for linear state-space systems can be applied to MIMO port-Hamiltonian systems in such a way as to preserve the \emph{port-Hamiltonian structure} in the reduced models, and preserve, as a consequence, \emph{passivity} as well. Moreover, we introduce a numerically efficient $\Htwo$-based algorithm for structure-preserving model reduction of port-Hamiltonian systems that produces high quality reduced order models in the general MIMO case. Numerical examples are presented to illustrate the effectiveness of the proposed method.

The rest of the paper is organized as follows:
 In Section \ref{sec:setup}, we review the solution to the rational tangential interpolation problem for general linear MIMO systems. Brief theory on port-Hamiltonian systems is given in Section \ref{sec:lphs}. Structure-preserving interpolatory model reduction of \phss in different coordinates is considered in Section \ref{sec:reduction_lphs} where we show theoretical equivalence of interpolatory reduction methods using different coordinate representations of port-Hamiltonian systems and discuss which is most robust and numerically effective. 
$\htwo$-based model reduction for port-Hamiltonian systems together with the proposed algorithm is presented in Section \ref{optimalh2ph} followed by numerical examples in Section \ref{sec:examples}.

\section{Interpolatory Model Reduction}
\label{sec:setup}
Let $\bfG(s)$ be a dynamical system with a state-space realization given as
\begin{eqnarray}
\label{eq:linsys}
\bfG(s) &:& \qquad \left\{\begin{array}{rcl}
\bfE \dot \bfx(t)& = & \bfA\bfx(t) + \bfB\, \bfu(t), \\
\bfy(t)& = & \bfC\, \bfx(t),
\end{array} \right.
\end{eqnarray}
where $\bfA,\bfE \in \IR^{n \times n}$,  $\bfE$ is nonsingular, 
$\bfB \in \IR^{n \times m}$, and $\bfC \in \IR^{p \times n}$.
In (\ref{eq:linsys}), for each $t$ we have: $\bfx(t)\in \IR^n$, $\bfu(t)\in \IR^m$  and $\bfy(t)\in \IR^p$ 
denoting, respectively, the \emph{state},  the \emph{input},
and  the \emph{output} of $\bfG(s)$. By taking the Laplace transform of (\ref{eq:linsys}),
we obtain the associated transfer function $\bfG(s)  =  \bfC (s\bfE-\bfA)^{-1}\bfB$.
Following  usual conventions,
the underlying dynamical system and its transfer function will both be denoted by
$\bfG(s)$.


We wish to produce a surrogate dynamical system $\bfG_r(s)$ of much smaller order  with a similar  state-space form
\begin{eqnarray}
\bfG_r(s) &:& \left\{ \begin{array}{l}   \bfE_r    \dot{\bfx}_r (t)
=  \bfA_r \bfx_r (t) + \bfB_r u(t) \\
        \bfy_r(t)  =  \bfC_r \bfx_r (t),  \end{array}
                   \right. \label{redsysintro}
\end{eqnarray}
where  $\bfA_r, \bfE_r \in \IR^{r \times r}$, $\bfB_r \in \IR^{r \times m}$, and $\bfC_r \in \IR^{p \times r}$  with $r \ll n$ such that the reduced-order system output
$\bfy_r(t)$ approximates the original system output, $\bfy(t)$, with high fidelity, as measured in an appropriate sense.
The appropriate sense here will be with respect to the $\Htwo$ norm and will be discussed  in Sections \ref{optimalh2} and \ref{optimalh2ph}.

We construct reduced order
models $\bfG_r(s)$ via Petrov-Galerkin projections. We choose
two $r$-dimensional  subspaces  $\mathcal{V}_r$ and   $\mathcal{W}_r$ and associated basis matrices
$\bfV_r \in \IR^{n \times r}$ and $\bfW_r \in \IR^{n \times r}$  (such that $\mathcal{V}_r = \rm{Range}(\bfV_r)$ and $\mathcal{W}_r = \rm{Range}(\bfW_r)$).  System dynamics are approximated as follows: we approximate
the full-order state $\bfx(t)$ as $\bfx(t) \approx \bfV_r \bfx_r(t)$ and force 
a Petrov-Galerkin condition as
\begin{eqnarray} \label{eqn:pg}
\mathbf{W}_{\! r}^{T}\left(\bfE \bfV_{\! r}\dot{\bfx}_r-\bfA\bfV_{\! r}\bfx_r -\bfB\,\bfu\right)=\mathbf{0},
\quad \bfy_r = \bfC\bfV_r \bfx_r,
\end{eqnarray}
 leading to a reduced-order model as in (\ref{redsysintro}) with
{
\begin{equation}  \label{red_projection}
\begin{array}{c}
\bfE_r = \bfW_r^T \bfE \bfV_r, \quad
\bfB_r = \bfW_r^T\bfB,  \\
\bfA_r = \bfW_r^T \bfA \bfV_r, \quad \mbox{and} \quad
\bfC_r = \bfC\bfV_r.
\end{array}
\end{equation}
}
The quality of the reduced model depends solely on the selection of the two subspaces
 $\mathcal{V}_r$ and   $\mathcal{W}_r$ (or equivalently on the choices of $\bfV_r$ and $\bfW_r$). 
 We choose $\bfV_r$ and $\bfW_r$ to force interpolation.
\subsection{Interpolatory projections}
The construction of reduced order models with interpolatory projections was initially
introduced by Skelton
{\it et. al.} in \cite{devillemagne1987mru,yousouff1984cer,yousuff1985lsa} and
later put into a robust numerical framework by
Grimme \cite{grimme1997krylov}. This
problem framework 
was recently adapted to MIMO dynamical systems of the form (\ref{eq:linsys}) by Gallivan {\it et al.} \cite{gallivan2005mrm}.

Beattie and Gugercin  in \cite{beattie2009ipm}
 later extended this framework to a much larger class of transfer functions, namely  those having a generalized coprime factorization
of the form $\bfH(s) = \bfC(s)\bfK(s)^{-1}\bfB(s) $ with $\bfB(s)$, $ \bfC(s)$, and $\bfK(s)$  given as meromorphic matrix-valued functions.
Given a full-order system (\ref{eq:linsys}) and interpolation points
$\{s_k\}_{k=1}^r \in \IC$ with corresponding (right) tangent directions 
$\{\bfsfb_k\}_{k=1}^r \in \IC^m$, the tangential interpolation problem
seeks a reduced-order system $\bfG_r(s)$ that interpolates $\bfG(s)$ at the selected points 
along the selected directions; i.e.,
\begin{equation} \label{eqn:rt}
\bfG(s_i) \bfsfb_i =  \bfG_r(s_i) \bfsfb_i,~~~
\hspace{0.5cm}{\rm for}~~i=1,\cdots,r, 
\end{equation}
Analogous left tangential interpolation conditions may be considered, however (\ref{eqn:rt}) will suffice for our purposes. 
Conditions forcing (\ref{eqn:rt}) to be satisfied by a reduced system of the form 
 (\ref{red_projection}) are provided by the following theorem (see \cite{gallivan2005mrm}).
\begin{thm}
\label{thm:interpolation}
Suppose $\bfG(s) = \bfC (s\bfE - \bfA)^{-1}\bfB$.  Given a set of distinct interpolation points $\{s_k\}_{k=1}^r$ and right tangent directions $\{\bfsfb_k\}_{k=1}^r$, define $\bfV_{\! r} \in \IC^{n \times r}$ as
\begin{eqnarray}
\bfV_{\! r} =
\left[(s_1\bfE-\bfA)^{-1}\bfB\bfsfb_1,~\cdots,~(s_r\bfE-\bfA)^{-1}\bfB\bfsfb_r\right] \label{eqn:Vr}
\end{eqnarray}
Then for any $\bfW_r\in \IC^{n\times r}$, the reduced order system $\bfG_r(s)= \bfC_r (s\bfE_r - \bfA_r)^{-1}\bfB_r$ defined as in (\ref{red_projection}) satisfies (\ref{eqn:rt})
 provided that
 $s_i\bfE-\bfA$ and  $s_i\bfE_r-\bfA_r$ are invertible for $i=1,\ldots,r$.
 \end{thm}

\begin{rem} \label{mimo_vs_siso}

Attempting interpolation of the full transfer function matrix (in all input and output directions) generally inflates the reduced system order by a factor of $m$ (the input dimension).  However, this is neither necessary nor desirable.  Effective reduced order models, indeed, $\Htwo$-\emph{optimal} reduced order models,  can be generated by forcing interpolation along selected tangent directions. 
  \end{rem}
 \begin{rem} \label{rem:Vderiv}
 The result of Theorem \ref{thm:interpolation} generalizes to higher-order interpolation (analogous to generalized Hermite interpolation) as follows:
 For a point $\hat{s}\in\IC$ and tangent direction $\hat{\bfsfb}$, suppose
 \vspace{-.2in}
{\small
\begin{equation}  \label{eqn:Vhigher}
\begin{array}{c}
\left(\left(\hat{s} \bfE - \bfA\right)^{-1} \bfE\right)^{k-1} (\hat{s} \bfE - \bfA)^{-1} \bfB \hat{\bfsfb} \in {\rm Range}(\bfV_r)\\
\mbox{for }k = 1, \ldots, N. 
\end{array}
\end{equation}
}
Then, for any $\bfW_r\in \IC^{n\times r}$, the reduced order system $\bfG_r(s)$ defined as in (\ref{red_projection}) satisfies
\begin{equation}  \label{eq:Hhigher}
\bfG^{(\ell)}(\hat{s}) \hat{\bfsfb}  =  \bfG^{(\ell)}_r(\hat{s}) \hat{\bfsfb},~\mbox{for }~\ell=0,\cdots,N-1,
\end{equation}
 provided that
 $\hat{s}\bfE-\bfA$ and  $\hat{s}\bfE_r-\bfA_r$ are invertible where $\bfG^{(\ell)}(s)$ denotes the $\ell^{\rm th}$ derivative of $\bfG(s)$. By combining this result with Theorem \ref{thm:interpolation}, one can match different interpolation orders for different interpolation points (analogous to generalized Hermite interpolation). For details; see, for example, \cite{gallivan2005mrm,antoulas2010imr}.
  \end{rem}
\begin{rem} \label{rem:basis}
Notice that in Theorem \ref{thm:interpolation} and Remark \ref{rem:Vderiv}, what guarantees the interpolation
is the range of the matrix $\bfV_r$, not the specific choice of basis in (\ref{eqn:Vr}). Hence,
one may substitute $\bfV_r$  with any matrix $\widehat{\bfV}_r$ having the same range.
In other words, for any $\widehat{\bfV}_r$ satisfying $\widehat{\bfV}_r = \bfV_r \bfL$ where $\bfL \in \IR^{r \times r}$ invertible, the interpolation property still holds true. This is a simple
consequence of  the fact that the basis change $\bfL$  from $\bfV_r$ to $\widehat{\bfV}_r$ amounts to a state-space transformation in the reduced model. 
We will make use of this property 
in discussing the effect of different coordinate representations for port-Hamiltonian systems. 
When the interpolation points occur in complex conjugate pairs (as always occurs in practice), then 
the columns of $\bfV_r$ also occur in complex conjugate pairs  and may be replaced with a \emph{real} basis,
instead.   We write $\bfV_r=[\![\bfv_1,\,\ldots,\,\bfv_r]\!]$ to represent that a real basis for $\bfV_r$ is chosen. 
Notice then that  $\hbfVr$ will also be a real matrix which then leads to real reduced order quantities as well.
\end{rem}


\subsection{Interpolatory optimal $\Htwo$ approximation}
\label{optimalh2}
 The selection of $r$ interpolation points $\{s_i\}_{i=1}^r \in \IC$ and corresponding tangent directions 
$\{\bfsfb_i\}_{i=1}^r \in \IC^m$ completely determines projecting subspaces and so is the sole factor determining the quality of a reduced-order models. Until recently  only heuristics were available to guide this process and the lack of a systematic selection process for interpolation points and tangent directions had been a key drawback to interpolatory model reduction. However,
 Gugercin {\it et  al.} \cite{gugercin2008hmr}
introduced a framework for determining optimal interpolation points and tangent directions  to find
optimal $\htwo$ reduced-order approximations, and also proposed an algorithm, the Iterative
Rational Krylov Algorithm (\emph{IRKA}) to compute the associated reduced-order models.  In this section we briefly review the interpolatory optimal $\htwo$ problem for the general linear systems 
$\bfG(s) = \bfC ( s\bfE - \bfA)^{-1}\bfB$ without structure.

Let $\mathcal{M}(r)$ denote the set of reduced order models having state-space dimension $r$ as in (\ref{redsysintro}).   Given a full order system $\bfG(s) = \bfC (s\bfE - \bfA)^{-1}\bfB$, the optimal-$\htwo$ reduced-order approximation to  $\bfG(s)$ of order $r$ is a solution to 
\begin{equation} \label{h2opt}
 \left\| \bfG - \bfG_r \right\|_{\htwo} =  \min_{ \mbox{\tiny{$
\widetilde{\bfG}_{r}\!\in\!\mathcal{M}(r)$}} } \left\| \bfG - \widetilde{\bfG}_{r} \right\|_{\htwo}
\end{equation}
where
$$
\left\| \bfG \right\|_{\htwo} :=
\left(\frac{1}{2\pi}\int_{-\infty}^{\infty}
  \left\| \bfG(\imath \omega )  \right\|_{\rm F}^2 d\omega\right)^{1/2}.
$$

There are a variety of methods that have been developed to solve (\ref{h2opt}).  They can be categorized as Lyapunov-based  methods such as
\cite{yan1999anapproximate,spanos1992anewalgorithm,halevi1992fwm,hyland1985theoptimal,wilson1970optimum,zigic1993contragredient}
and interpolation-based optimal $\htwo$ methods such as
\cite{meieriii1967approximation,gugercin2008hmr,gugercin2006rki,vandooren2008hom,bunse-gerstner2009hom,gugercin2005irk,kubalinska2007h0i,beattie2007kbm,beattie2009trm}.
Although both frameworks are theoretically equivalent \cite{gugercin2008hmr}, 
interpolation-based optimal $\htwo$ methods carry important computational advantages.
Our focus here will be  on interpolation-based approaches. 

The optimization problem (\ref{h2opt}) is nonconvex, so obtaining a global minimizer is difficult in the best of 
circumstances.  Typically local minimizers are sought instead and as a practical matter, 
 the usual approach is to find reduced-order models satisfying the first-order necessary conditions of 
 $\htwo$-optimality.
 
Interpolation-based $\htwo$-optimality conditions for SISO systems
were introduced by
Meier and Luenberger \cite{meieriii1967approximation}.
Interpolatory $\htwo$-optimality conditions for MIMO systems
 have recently been developed
by \cite{gugercin2008hmr,bunse-gerstner2009hom,vandooren2008hom}.
\begin{thm}  \label{h2cond}  \label{meier}
Given  a full-order system $\bfG(s) =  \bfC (s \bfE - \bfA)^{-1}\bfB$, let $\bfG_r(s)=\sum_{i=1}^r\frac{1}{s-\hat{\lambda}_i}\bfsfc_i\bfsfb_i^T$
 be  the best $r^{\rm th}$ order
 approximation of $\bfG$ with respect to the $\htwo$ norm.
Then
\begin{align}     \label{H2optcond1}
  \mbox{(a)}\ & \bfG(-\hat{\lambda}_k) \bfsfb_k =\bfG_r(-\hat{\lambda}_k) \bfsfb_k,
\\ \mbox{(b)}\ &  \bfsfc_k^T \bfG(-\hat{\lambda}_k) =  \bfsfc_k^T\bfG_r(-\hat{\lambda}_k), \mbox{ and }\quad
   \label{H2optcond2}\\
  \mbox{(c)}\   &\bfsfc_k^T \bfG'(-\hat{\lambda}_k) \bfsfb_k
   =  \bfsfc_k^T\bfG_r'(-\hat{\lambda}_k) \bfsfb_k
   \label{H2optcond3}
\end{align}
for $k=1,\,2,\,...,\,r.$
\end{thm}
Theorem \ref{h2cond} asserts  that any solution to the optimal-$\htwo$ problem (\ref{h2opt})
 must be a bitangential Hermite interpolant to $\bfG(s)$ at the mirror images of the reduced system poles. 
This would be a straightforward construction, if the reduced system poles and residues were known \emph{a priori}.
However, they are not and the \emph{Iterative Rational Krylov Algorithm}  \emph{IRKA}  \cite{gugercin2008hmr} resolves the problem
by iteratively correcting the interpolation points by reflecting the
 current reduced-order system poles across the imaginary axis to become the next set of interpolation points  
 and the directions.  The next tangent directions are residue directions taken from the current reduced model. Upon convergence, the necessary conditions of Theorem  \ref{h2cond} are met.
\emph{IRKA} corresponds to applying interpolatory model reduction iteratively using   
{\small
\begin{equation} \label{eqn:VrWr}
\begin{array}{l}
\bfV_{\! r} = 
\left[(s_1\bfE-\bfA)^{-1}\bfB\bfsfb_1,\cdots,(s_r\bfE-\bfA)^{-1}\bfB\bfsfb_r\right] \\ 
\bfW_{\! r} =
\left[(s_1\bfE-\bfA)^{-T}\bfC^T\bfsfc_1,\cdots,(s_r\bfE-\bfA)^{-T}\bfC^T\bfsfc_r\right]
\end{array}
\end{equation}
}
where $s_i$, $\bfsfb_i$ and $\bfsfc_i$ are, respectively, the interpolation points, and right and left tangent directions at step $k$ of \emph{IRKA}.   Note that the Hermite bitangential interpolation conditions 
(\ref{H2optcond1})-(\ref{H2optcond3}) enforce a choice on $\bfW_r$ as in (\ref{eqn:VrWr}) in contrast to Theorem \ref{thm:interpolation} where $\bfW_r$ could be chosen arbitrarily since bitangential interpolation is not required. We revisit this issue in Section~\ref{optimalh2ph}.
For details on  \emph{IRKA}, see \cite{gugercin2008hmr}.


\section{Linear port-Hamiltonian Systems}
\label{sec:lphs}

In the absence of algebraic constraints, linear port-Hamiltonian systems take the following form (\cite{Schaft00,geoplexbook2009,PolvdSc09})
\begin{eqnarray}  \label{eq:lphs}
\bfG(s) &:&\qquad \
\left\{\begin{array}{rcl}
\dot \bfx &=&  (\bfJ - \bfR)\bfQ\bfx + \bfB\,\bfu, \\[1mm]
\bfy &=&  \bfB^T\bfQ\bfx,
\end{array}\right.
\end{eqnarray}
where $\bfJ = -\bfJ^T$, $\bfQ = \bfQ^T$, and $\bfR = \bfR^T\geqslant 0$.   $\bfQ$ is the \emph{energy} matrix;
$\bfR$  is the \emph{dissipation} matrix; $\bfJ$ and $\bfB$ determine the \emph{interconnection} structure of the system.    $H(\bfx) = \frac{1}{2}\bfx^T\bfQ\bfx$ defines the \emph{total energy} (the \emph{Hamiltonian}). For all cases of interest, $H(\bfx)\geq 0$.
Note that the system \eqref{eq:lphs} has the form \eqref{eq:linsys} with
$\bfE = \bfI$, $\bfA = (\bfJ -\bfR)\bfQ$, and $\bfC = \bfB^T\bfQ$.

The state variables $\bfx \in \mathbb{R}^n$ are called \emph{energy} 
variables, since the total energy $H(\bfx)$ is expressed as a function of these variables.   
 $\bfu, \bfy \in \mathbb{R}^m$ are called \emph{power} variables, since the scalar product $\bfu^T \bfy$ equals the power supplied to the system.
Note that 
\begin{equation*} 
\frac{d}{dt}\ \mbox{{$\frac{1}{2}$}}\,\bfx^T\bfQ\bfx\ = \ \bfu^T \bfy - \bfx^T\bfQ\bfR\bfQ\bfx\ \leqslant\ \bfu^T \bfy
\end{equation*}
--- that is, port-Hamiltonian systems are \emph{passive}  and the change in total system energy is bounded by the work done on the system. 

We will assume throughout that $\bfQ$ is positive definite.  In the case of zero dissipation ($\bfR = \mathbf{0}$), the full-order model $\bfG(s)$ of \eqref{eq:lphs} is both lossless and stable.  For positive semidefinite $\bfR$, $\bfG(s)$ could be either stable or asymptotically stable in general; and for $\bfR$ strictly positive definite, the system is asymptotically stable. 
We will assume that $\bfG(s)$ is merely stable except in Section \ref{optimalh2ph} where  we assume asymptotic stability in order 
to assure a bounded $\Htwo$ norm.

%
%
%
%
%
%


\begin{example}
\label{exm:ladder_x}
Consider the Ladder Network illustrated in Fig. 1, with $C_1, C_2, L_1, L_2, R_1, R_2, R_3$ being, respectively, the capacitances,
inductances, and resistances of idealized linear circuit elements described in the figure.
The port-Hamiltonian representation of this physical system has the form \eqref{eq:lphs} with
\begin{eqnarray}
\label{eq:ladder_ABC}
\begin{array}{rcl}
\bfQ &=& \mathsf{diag}(C_1^{-1}, L_1^{-1},C_2^{-1}, L_2^{-1}), \\ \\
\bfJ &=&\mathsf{tridiag}\!\left(
\begin{array}{ccccccc}
& -1 && -1 && 1 & \\[-.1in]
0 && 0 && 0 && 0 \\[-.1in]
& 1 && 1 && -1 &
\end{array}
 \right), \\
\bfR &=& \mathsf{diag}(0, R_1, 0, R_2+R_3),
\quad \mbox{and}\quad \bfB = \left[\begin{array}{cc} 1 & 0 \\ 0 & 0 \\  0 & 0 \\ 0 & 1  \end{array}\right].
\end{array}
\end{eqnarray}
The $\bfA$ matrix is
\begin{eqnarray*}
\bfA &=&
\left[\begin{array}{cccc}
0 & -L_1^{-1} & 0 & 0\\
C_1^{-1} & -R_1L_1^{-1} & -C_2^{-1} & 0 \\
0 & L_1^{-1} & 0 & L_2^{-1} \\
0 & 0 & -C_2^{-1} & -(R_2+R_3)L_2^{-1}
\end{array}\right] \\
\end{eqnarray*}
and the state-space vector $\bfx$ is given as
\begin{eqnarray*}
\bfx^T = \left[\begin{array}{cccc} q_1 & \phi_1 & q_2 & \phi_2 \end{array}\right]
\end{eqnarray*}
with $q_1, q_2$ being the charges on the capacitors $C_1, C_2$ and $\phi_1, \phi_2$ being
the fluxes over the inductors $L_1, L_2$ correspondingly. The inputs of the
system, $\{u_1, u_2\}$ are given by the current $I$ on the left side and the voltage $U$ on the right side
of the Ladder Network. The \ph outputs $\{y_1, y_2\}$ are
the voltage over the first capacitor $U_{C_1}$ and the
current through the second inductor $I_{L_2}$.

The system matrices $\bfA$, $\bfB$, $\bfC$, and $\bfE$ of (\ref{eq:linsys}) follow directly from writing the linear input-state differential equation for this system. $\bfQ$ may be derived from the Hamiltonian $H(\bfx) = \frac{1}{2}\bfx^T\bfQ\bfx$. Once $\bfA$ and $\bfQ$ are known, it is easy to derive the dissipation matrix, $\bfR$,
 and the structure matrix, $\bfJ$, corresponding to the Kirchhoff laws, 
such that $\bfA = (\bfJ - \bfR)\bfQ$. The \ph output matrix $\bfC$ is given as
\begin{equation}
\label{eq:ladder_ABCC}
\bfC =\bfB^T\bfQ=\left[ \begin{array}{cccccc}
   \frac{1}{C_1} & 0 & 0 & 0 \\
   0 & 0 & 0 &    \frac{1}{L_2} \\
\end{array}  \right].
\end{equation}
\end{example}


\begin{figure}[t]
\label{fig:4ladnet}
\centering
\includegraphics[width=8.5cm]{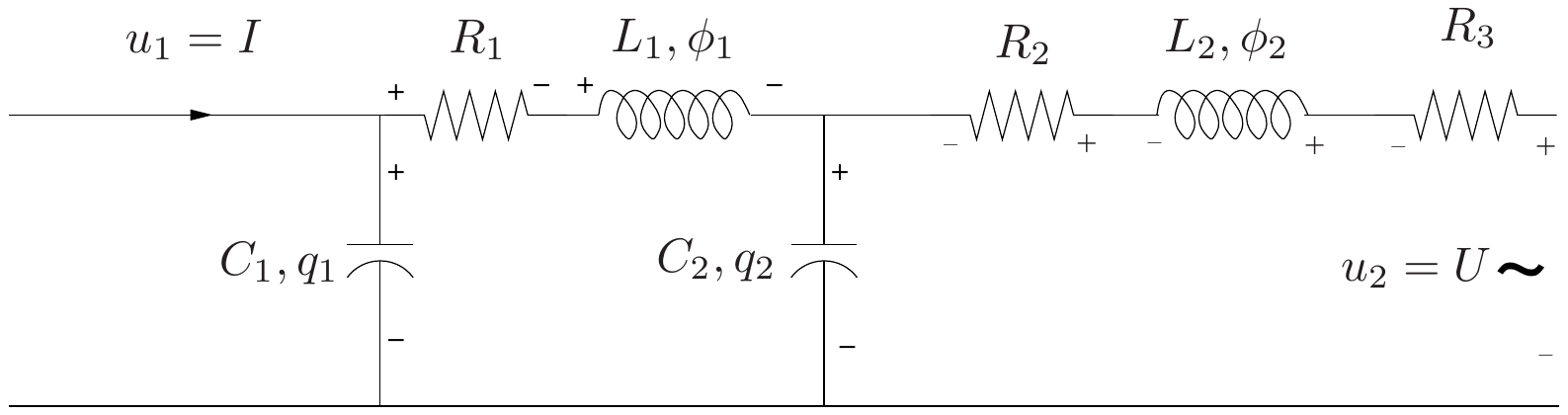}
\caption{Ladder Network}
\end{figure}

  Recall from \cite{geoplexbook2009,Polyuga_2010_thesis,PolvdSc09} that the port-Hamiltonian system \eqref{eq:lphs} may be represented in  \emph{co-energy coordinates} as
\begin{eqnarray} 
\label{eq:lphs_coen}
\left\{\begin{array}{rcl}
\dot \bfe & = & \bfQ(\bfJ - \bfR)\bfe + \bfQ\bfB\,\bfu,\\
\bfy & = & \bfB^T\bfe.
\end{array}\right.
\end{eqnarray}
The coordinate transformation \cite{geoplexbook2009,Polyuga_2010_thesis} between energy coordinates, $\bfx$, and co-energy coordinates, $\bfe$, is given by
\begin{equation} \label{coenergyCoord}
\bfe = \bfQ\bfx.
\end{equation}


\begin{example}
\label{exm:ladder_e}
(continued). The co-energy state vector, $\bfe$, for the Ladder Network from Example \ref{exm:ladder_x} is given as
\begin{eqnarray*}
\bfe^T = \left[\begin{array}{cccc} U_{C_1} & I_{L_1} & U_{C_2} & I_{L_2} \end{array}\right]
\end{eqnarray*}
with $U_{C_1}, U_{C_2}$ being the voltages on the capacitors $C_1, C_2$ and $I_{L_1}, I_{L_2}$
being the currents through the inductors $L_1, L_2$, respectively.
\end{example}

\subsection{Balancing for port-Hamiltonian systems}
\label{sec:effort}

To make the presentation self-contained, we review a recent structure-preserving, balancing-based model reduction method for port-Hamiltonian systems, the \emph{effort-constraint method} (\cite{Polyuga_vdSchaft2011_Effort_flow,vdSchaft_Polyuga_09cdc,Polyuga_2010_thesis,PolvdSc08,PolvdSc_2011_lecture_notes}), 
which will be used for comparisons in our  numerical examples.

Consider a full order port-Hamiltonian system realized as in \eqref{eq:lphs} with respect to energy coordinates.    
Consider the associated \emph{balancing transformation}, $\bfT_b$, defined in the usual way
(see \cite{antoulas2005approximation} for a complete treatment): $\bfT_b$ simultaneously diagonalizes the observability Gramian,
$\mathbb{G}_o$, and  the controllability Gramian,
$\mathbb{G}_c$, so that
$$
\bfT_b^{-1} \mathbb{G}_o \bfT_b^{-T} = \bfT_b^{T} \mathbb{G}_c \bfT_b = \mathsf{diag}(\sigma_1,\,\sigma_2,\,
\ldots,\,\sigma_n).
$$
where $\sigma_1,\,\sigma_2,\,
\ldots,\,\sigma_n$ are the Hankel singular values of the system $\bfG(s)$.  This balancing transformation is the same that is used in the context of regular \emph{balanced truncation}.

The state-space transformation associated with balancing (as with any other linear coordinate transformation), will preserve 
port-Hamiltonian structure of the system \eqref{eq:lphs}.
Indeed, if we define \emph{balanced coordinates},  as $\bfx_b = \bfT_b\bfx$, then we have
\begin{eqnarray}
\label{eq:lphs_balanced}
\left\{\begin{array}{lcl}
\dot \bfx_b & =& (\bfJ_b - \bfR_b)\bfQ_b\bfx_b + \bfB_b\bfu, \\[1mm]
\bfy & =& (\bfB_b)^T\bfQ_b\bfx_b,
\end{array}\right.
\end{eqnarray}
where $\bfT_b\bfR \bfT_b^{T} = \bfR_b = \bfR_b^T\geq 0$ is the dissipation matrix, $\bfT_b\bfJ \bfT_b^{T} = \bfJ_b = -\bfJ_b^T$ is the structure matrix and $\bfT_b^{-T}\bfQ \bfT_b^{-1} = \bfQ_b = \bfQ_b^T> 0$ is the energy matrix in balanced coordinates. In this case, $\bfB_b = \bfT_b\bfB$.

Now split the state-space vector, $\bfx_b$, into dominant and codominant components:
$\displaystyle \bfx_b=\begin{bmatrix}
\bfx_{b1} \\
\bfx_{b2}
\end{bmatrix}, $
where $\bfx_{b1} \in \IR^r$. Regular balanced truncation proceeds by truncating the codominant state-space components, in effect forcing the system to evolve under the constraint $\bfx_{b2}=\mathbf{0}$.  This will destroy port-Hamiltonian structure in the remaining reduced order system that determines the evolution of $\bfx_{b1}$.

%

In contrast to this, the \emph{effort-constraint method} produces the following reduced order port-Hamiltonian system
\begin{equation} \label{eq:RedSysEfConstr}
\begin{array}{ccl}
\dot \bfx_{b1}  &=&  (\bfJ_{b11}  -  \bfR_{b11})  \bfS_b \bfx_{b1}  + \bfB_{b1}\bfu,\\[1mm]
\bfy_r  &=&  (\bfB_{b1})^T \bfS_b\bfx_{b1},
\end{array}
\end{equation}
where $\bfS_b = \bfQ_{b11}-\bfQ_{b12}(\bfQ_{b22})^{-1}\bfQ_{b21}$ is the Schur complement of the energy matrix $\bfQ_b$ in balanced coordinates; the other matrices are of corresponding dimension.  The relationship of the effort-constraint method to the reduction of the underlying (full-order)  Dirac structure is discussed in \cite{Polyuga_vdSchaft2011_Effort_flow,vdSchaft_Polyuga_09cdc,Polyuga_2010_thesis}. A more direct approach is given in \cite{PolvdSc_2011_lecture_notes};  another approach using scattering coordinates can be found in \cite{PolvdSc08}.

\begin{rem}
\label{rem:ec_vs_truncation}
The effort-constraint method can be formulated with Petrov-Galerkin projections \cite{Polyuga_vdSchaft2011_Effort_flow,vdSchaft_Polyuga_09cdc,Polyuga_2010_thesis}.
Even though a balancing transformation is used in developing the effort-constraint method as expressed in \eqref{eq:RedSysEfConstr}, the method is not equivalent to the usual balanced truncation method. Note that balanced truncation does not preserve \ph structure.
For more details, see \cite{Polyuga_2010_thesis}.
\end{rem}


\section{Interpolatory model reduction of port-Hamiltonian systems}
\label{sec:reduction_lphs}

Now we consider how to create an interpolatory reduced-order model that retains port-Hamiltonian structure inherited from the original system, thus guaranteeing  that it is both stable and passive. 

Note that Theorem \ref{thm:interpolation} may be applied to a port-Hamiltonian system with an arbitrary choice of $\bfW_r$. As a result, $\bfA_r$ will have the form $\bfW_r^T(\bfJ-\bfR)\bfQ\bfV_r$, which can be no longer represented as $(\bfJ_r-\bfR_r)\bfQ_r$ with a skew-symmetric $\bfJ_r$, and
symmetric positive semi-definite $\bfR_r$ and
positive definite
$\bfQ_r$ matrices. Below, we  develop several approaches that will choose $\bfW_r$ carefully to allow the structure preservation.
\subsection{Interpolatory projection with respect to energy coordinates} \label{sec:rom_orig}
We first show how to achieve interpolation and preservation of port-Hamiltonian structure simultaneously in the original energy coordinate representation. The choice of $\bfW_r$ plays a crucial role.
\begin{thm}  \label{thm:origcoor}
Suppose $\bfG(s)$ is a linear port-Hamiltonian system, as  in \eqref{eq:lphs}.
Let $\{s_i\}_{i=1}^r \subset \IC$ be a set of $r$ distinct interpolation points
with corresponding tangent directions $\{\bfsfb_i\}_{i=1}^r \subset \IC^m$, both sets being closed under conjugation.
Construct $\bfV_r$
as in (\ref{eqn:Vr}) using
$\bfA  = (\bfJ - \bfR)\bfQ$ and $\bfE = \bfI$, 
so that\\[-.2in]
{\small 
$
\bfV_r = [\![(s_1 \bfI - (\bfJ-\bfR)\bfQ)^{-1}\bfB\bfsfb_1,\ldots,(s_r \bfI - (\bfJ-\bfR)\bfQ)^{-1}\bfB\bfsfb_r]\!].
$
}
 
Let $\bfM\in\IC^{r\times r}$ be any nonsingular matrix such that $\hbfVr = \bfVr\bfM$ is real.   

Define $\hbfWr=\bfQ\hbfVr(\hbfVr^T\bfQ \hbfVr)^{-1}$ and then
\begin{equation} \label{eqn:lphs_rom_ss}
\begin{array}{ll}
\bfJ_r = \hbfWr^T\bfJ\hbfWr,  &\quad  \bfQ_r =  \hbfVr^T\bfQ \hbfVr, \\
\bfR_r = \hbfWr^T\bfR\hbfWr, &\quad \mbox{and}\quad \bfB_r = \hbfWr^T\bfB. \end{array}
\end{equation}
%
Then, the reduced-order model
\begin{equation} 
\label{eq:lphs_rom}
\bfG_r(s): \quad \left\{
\begin{array}{c}
 \dot \bfx_r  =  (\bfJ_r - \bfR_r)\bfQ_r\bfx_r + \bfB_r\,\bfu\\
\bfy_r  =  \bfB_r^T \bfQ_r\bfx_r
\end{array}
\right.
\end{equation}
is port-Hamiltonian, passive, and 
$$
\bfG(s_i) \bfsfb_i =  \bfG_r(s_i) \bfsfb_i,\qquad \mbox{ for }i=1,\cdots,r.
$$
That is, $\bfG_r(s)$ interpolates $\bfG(s)$ at $\{s_i\}_{i=1}^r$ along the tangent directions $\{\bfsfb_i\}_{i=1}^r$.
 \end{thm}
{\bf Proof:}  
In light of Remark \ref{rem:basis}, there is at least one (and hence, an infinite number) of possible 
$\bfM\in\IC^{r\times r}$ satisfying the conditions of the Theorem and so, $\hbfVr\in \IR^{n\times r}$. 
Trivially, one may observe that the reduced-order system  (\ref{eq:lphs_rom}) is real and retains port-Hamiltonian structure with a positive definite $\bfQ_r$, and so must be passive.    To show that 
$\bfG_r(s)$ interpolates $\bfG(s)$ 
note, as before, that the port-Hamiltonian system \eqref{eq:lphs} has the standard form \eqref{eq:linsys} with
$\bfE = \bfI$, $\bfA = (\bfJ -\bfR)\bfQ$, and $\bfC = \bfB^T\bfQ$. 
An interpolatory reduced model  may be defined as in (\ref{redsysintro}) using 
$\bfV_r\stackrel{\mbox{\tiny{\textsf{def}}}}{=}\hbfVr$ and 
$\bfW_r\stackrel{\mbox{\tiny{\textsf{def}}}}{=}\hbfWr$. We then have
$\bfE_r = \hbfWr^T \hbfVr=\bfI_r$ and
$\bfA_r = \hbfWr^T \bfA \hbfVr =  \hbfWr^T (\bfJ-\bfR)\bfQ \hbfVr$.   
From Theorem \ref{thm:interpolation}, this model  interpolates $\bfG(s)$ at $\{s_i\}_{i=1}^r$ along the tangent directions $\{\bfsfb_i\}_{i=1}^r$ as required but it is not obvious that this is the same reduced system $\bfG_r(s)$ that we have defined above. 

Note that $\hbfWr\hbfVr^T$ is a (skew) projection onto $\mbox{Range}(\bfQ\hbfVr)$ and so $\bfQ\hbfVr=\hbfWr\hbfVr^T\bfQ\hbfVr$. 
Thus, we have
\begin{align*}
   \bfA_r = &  \hbfWr^T (\bfJ-\bfR)\bfQ \hbfVr\\
       = &  \hbfWr^T (\bfJ-\bfR)\hbfWr\hbfVr^T\bfQ \hbfVr \\
       =&  (\hbfWr^T\bfJ \hbfWr-\hbfWr^T\bfR\hbfWr)\hbfVr^T\bfQ \hbfVr \\
      =&  (\bfJ_r-\bfR_r)\bfQ_r
\end{align*}
          and  
\begin{align*}
\bfC_r=&\bfC \hbfVr=\bfB^T\bfQ \hbfVr\\
=&\bfB^T\hbfWr\hbfVr^T\bfQ \hbfVr =\bfB_r^T\bfQ_r,
\end{align*}
which verifies that the port-Hamiltonian reduced system, $\bfG_r(s)$, interpolates 
$\bfG(s)$ at $\{s_i\}_{i=1}^r$ along $\{\bfsfb_i\}_{i=1}^r$ as required.
$\Box$.

\begin{rem}
Theorem \ref{thm:origcoor} can be generalized to include higher-order (generalized Hermite) interpolation following ideas described in Remark \ref{rem:Vderiv}.  Indeed,  given  an interpolation point $\hat{s}$ and tangent direction $\hat{\bfsfb}$,  augment $\bfV_r$ so that (\ref{eqn:Vhigher}) holds. The remaining part of the construction of Theorem  \ref{thm:origcoor} proceeds unchanged and the resulting reduced-order model will retain port-Hamiltonian structure (hence be both stable and passive) and will interpolate higher-order derivatives of $\bfG(s)$ as in (\ref{eq:Hhigher}).
\end{rem}

\subsection{The influence of state-space transformations}
Let $\bfT\in\IR^{n \times n}$ be an arbitrary invertible matrix representing a state-space transformation $\tilde{\bfx}=\bfT\bfx$.   As one may recall from the discussion of Section \ref{sec:effort}, port-Hamiltonian structure will be maintained under state-space transformations:
\begin{equation*} 
\begin{array}{rcl}
\dot{\bfx}& = & (\bfJ - \bfR)\bfQ \bfx + \bfB\,\bfu, \\
\bfy& = & \bfB^T \bfQ \bfx.
\end{array}
\iff
\begin{array}{rcl}
\dot{\tbfx}& = & (\tbfJ - \tbfR)\tbfQ\tbfx + \tbfB\,\bfu, \\
\bfy& = & \tbfB^T \tbfQ \tbfx.
\end{array}
\end{equation*}
with transformed quantities
\begin{equation*} 
\begin{array}{c}
\tbfJ = \bfT\bfJ \bfT^{T},~~
\tbfR  = \bfT \bfR \bfT^{T},\\
\tbfQ  = \bfT^{-T} \bfQ \bfT^{-1},\mbox{ and }
\tbfB =  \bfT \bfB.
\end{array}
\end{equation*}

\begin{thm}  \label{thm:equivalent}
Suppose $\bfG(s)$ is a port-Hamiltonian system as in \eqref{eq:lphs}, with two different port-Hamiltonian realizations connected via a state-space transformation, $\bfT$, as above.  
Given a set of $r$ distinct interpolation points $\{s_i\}_{i=1}^r \in \IC$ with corresponding tangent directions 
$\{\bfsfb_i\}_{i=1}^r \in \IC^m$ (closed under conjugation), then interpolatory reduced-order port-Hamiltonian models produced from either realization via Theorem \ref{thm:origcoor} will be identical to one another.
\end{thm}
{\bf Proof:} Define primitive interpolatory basis matrices for each realization as
{\small 
$
\bfV_r = [\![(s_1 \bfI - (\bfJ-\bfR)\bfQ)^{-1}\bfB\bfsfb_1,\ldots,(s_r \bfI - (\bfJ-\bfR)\bfQ)^{-1}\bfB\bfsfb_r]\!]
$
}

and
{\small 
$
\widetilde{\bfV}_r = [\![(s_1 \bfI - (\tbfJ-\tbfR)\tbfQ)^{-1}\tbfB\bfsfb_1,\ldots,(s_r \bfI - (\tbfJ-\tbfR)\tbfQ)^{-1}\tbfB\bfsfb_r]\!].
$
}
Elementary manipulations show that $\widetilde{\bfV}_r=\bfT\bfV_r $.  Note that any $\bfM\in\IC^{r\times r}$ for which 
$\widehat{\bfV}_r=\bfV_r\bfM$ is real makes $\widehat{\widetilde{\bfV}}_r=\widetilde{\bfV}_r\bfM$ real as well (and vice versa).  
Then we define, following the construction of Theorem \ref{thm:origcoor}, $\hbfWr=\bfQ\widehat{\bfV}_r$ and 
$\widehat{\widetilde{\bfW}}_r=\tbfQ \widehat{\widetilde{\bfV}}_r$.  Observe that $\widehat{\widetilde{\bfW}}_r= \bfT^{-T}\hbfWr$ and $\widehat{\widetilde{\bfV}}_r= \bfT\widehat{\bfV}_r$.  So,
\begin{equation*}
\begin{array}{c}
\bfJ_r = \hbfWr^T\bfJ\hbfWr = \widehat{\widetilde{\bfW}}_r^{\raisebox{-1.3ex}{\scriptsize{$T$}}}\tbfJ \widehat{\widetilde{\bfW}}_r= \widetilde{\bfJ}_r,\\
\bfQ_r =  \hbfVr^T\bfQ \hbfVr =  \widehat{\widetilde{\bfV}}_r^{\raisebox{-1.3ex}{\scriptsize{$T$}}}\tbfQ \widehat{\widetilde{\bfV}}_r= \widetilde{\bfQ}_r,\\
\bfR_r = \hbfWr^T\bfR\hbfWr= \widehat{\widetilde{\bfW}}_r^{\raisebox{-1.3ex}{\scriptsize{$T$}}}\tbfR\widehat{\widetilde{\bfW}}_r = \widetilde{\bfR}_r,\\
\quad \mbox{and}\quad \bfB_r = \hbfWr^T\bfB =\widehat{\widetilde{\bfW}}_r^{\raisebox{-1.3ex}{\scriptsize{$T$}}} \tbfB = \widetilde{\bfB}_r. \qquad \Box
\end{array}
\end{equation*}

\begin{rem} One may conclude from Theorem \ref{thm:equivalent} that prior to calculating a reduced-order interpolatory model, there may be little advantage in first applying a state-space transformation (e.g., a balancing transformation, $\bfT_b$, as described in Section \ref{sec:effort} or transforming to co-energy coordinates with $\bfT=\bfQ$ as in (\ref{coenergyCoord})).  Indeed, choosing a realization that exhibits advantageous sparsity patterns facilitating the linear system solves necessary to produce $\bfV_r$ will likely be the most effective choice.  If we choose  $\bfT$ to satisfy $\bfT^T\bfT=\bfQ$ (for example, if $\bfT$ is a Cholesky factor for $\bfQ$)  then, with respect to the new $\tilde{\bfx}$-coordinates (called ``scaled energy coordinates"), we find that $\tbfQ=\bfI$ and so 
$\widehat{\widetilde{\bfW}}_r= \widehat{\widetilde{\bfV}}_r$.  Notice that in this case, scaled energy coordinates and scaled co-energy coordinates become identical.  Notwithstanding these simplifications, unless  the original $\bfQ$ is diagonal (so that the transformation to scaled energy coordinates is cheap), there appears to be little justification for global coordinate transformations preceding the construction of an interpolatory reduced order model.  
\end{rem}
\section{$\htwo$-based reduction of port-Hamiltonian systems}
\label{optimalh2ph}
Inspired by the optimal $\Htwo$ model reduction method described in Section~\ref{optimalh2},
we propose here an algorithm that produces high quality  reduced-order port-Hamiltonian models of port-Hamiltonian systems. 
For the $\Htwo$ norm to be well defined, we assume that 
the full-order model $\bfG(s)$ is asymptotically stable.


Even though the optimality conditions (\ref{H2optcond1})-(\ref{H2optcond3}) can be satisfied for generic model reduction settings starting with a general system $\bfG(s) = \bfC(s\bfE-\bfA)^{-1}\bfB$ and reducing to $\bfG_r(s) = \bfC_r(s\bfE_r-\bfA_r)^{-1}\bfB_r$ without any restriction on structure,  it will generally not be possible to satisfy  all optimality conditions  while also retaining port-Hamiltonian structure; there are not enough degrees of freedom to achieve both goals simultaneously. 
Note that (\ref{H2optcond1})-(\ref{H2optcond3}) corresponds to $r(pm+1)$ conditions --- precisely the number of  degrees of freedom in 
$\bfG_r(s) = \bfC_r(s\bfE_r-\bfA_r)^{-1}\bfB_r  =\sum_{i=1}^r\frac{1}{s-\hat{\lambda}_i}\bfsfc_i\bfsfb_i^T.$
This representation is a general partial fraction expansion and need not correspond to a 
port-Hamiltonian system. Adding port-Hamiltonian structure on top of (\ref{H2optcond1})-(\ref{H2optcond3})  creates an overdetermined set of optimality conditions that cannot typically be satisfied.   Alternatively, note that the optimality conditions
(\ref{H2optcond1})-(\ref{H2optcond3})  requires choosing $\bfW_r$ in a specific way as in
(\ref{eqn:VrWr}) unlike the way it was chosen in Theorem \ref{thm:origcoor} in order to preserve structure.  
See also Remark \ref{PHandH2Opt}. 


We chose to enforce only a subset of first-order optimality conditions and use the remaining degrees of freedom to preserve structure. 
We  enforce (\ref{H2optcond1}) with a modifed version of \emph{IRKA} that maintains port-Hamiltonian structure at every step.   
A description of our proposed algorithm is given as Algorithm \ref{phirka}. 
\begin{alg} 
 {\bf IRKA for MIMO port-Hamiltonian systems (IRKA-PH):}
\begin{enumerate}
 \item  Let $\bfG(s) = \bfB^T\bfQ(s\bfI - (\bfJ-\bfR)\bfQ)^{-1}\bfB$ be as in \eqref{eq:lphs}.
\item Choose initial interpolation points $\{s_1,\ldots,s_r\}$ and tangential directions
$\{\bfsfb_1,\ldots, \bfsfb_r\}$; both sets closed under conjugation.
\item Construct\\
$\bfV_r = [\![(s_1 \bfI - (\bfJ-\bfR)\bfQ)^{-1}\bfB\bfsfb_1,\ldots,$\\
\hspace*{1.2in}$\ldots, (s_r \bfI - (\bfJ-\bfR)\bfQ)^{-1}\bfB\bfsfb_r]\!].$
\item Choose a nonsingular matrix $\bfM\in\IC^{r\times r}$  such that $\hbfVr = \bfVr\bfM$ is real. 
\item Calculate $\hbfWr=\bfQ\hbfVr(\hbfVr^T\bfQ \hbfVr)^{-1}$
\item until convergence
\begin{enumerate}
\item $\bfJ_r = \hbfWr^T\bfJ\hbfWr$, $\bfR_r = \hbfWr^T\bfR\hbfWr$, \\
$\bfQ_r =  \hbfVr^T\bfQ \hbfVr$ 
and $\bfB_r = \hbfWr^T \bfB$.
\item $\bfA_r = (\bfJ_r-\bfR_r)\bfQ_r$
  \item Compute
   $\bfA_r \bfx_i = \widehat{\lambda}_i\bfx_i,~~\bfy_i^*\bfA_r = \widehat{\lambda}_i\bfy_i^*$\\
  ${\rm with}~~\bfy_i^*\bfx_j = \delta_{ij}$
        for left and right eigenvectors $\bfy_i^*$ and $\bfx_i$  associated with  $\widehat{\lambda}_i$.
\item $s_i \leftarrow -\widehat{\lambda}_i$ and $\bfsfb_i^T \leftarrow \bfy_i^*\bfB_r$ for $i=1,\ldots,r$.
\item Compute 
$\bfV_r = [\![(s_1 \bfI - (\bfJ-\bfR)\bfQ)^{-1}\bfB\bfsfb_1,\ldots,$\\
\hspace*{1.2in}$\ldots, (s_r \bfI - (\bfJ-\bfR)\bfQ)^{-1}\bfB\bfsfb_r]\!].$
\item Choose a nonsingular matrix $\bfM\in\IC^{r\times r}$  such that $\hbfVr = \bfVr\bfM$ is real. 
\item Calculate $\hbfWr=\bfQ\hbfVr(\hbfVr^T\bfQ \hbfVr)^{-1}$
\end{enumerate}
\item The final reduced model is given by  \\
$\bfJ_r = \hbfWr^T\bfJ\hbfWr$, $\bfR_r = \hbfWr^T\bfR\hbfWr$,
 $\bfB_r =\hbfWr^T\bfB$, \\
\hspace*{.4in} $\bfQ_r = \hbfVr^T\bfQ \hbfVr$, and $\bfC_r = \bfB_r^T \bfQ_r$.
\end{enumerate}
\label{phirka}
\end{alg} 
\begin{thm}
Let $\bfG(s) = \bfB^T\bfQ(s\bfI - (\bfJ-\bfR)\bfQ)^{-1}\bfB$ be a port-Hamiltonian system as in \eqref{eq:lphs}.
Suppose \emph{IRKA-PH} as described in Algorithm \ref{phirka} converges to  a reduced model
\begin{equation} \label{partFracExp}
\bfG_r(s)=\sum_{i=1}^r\frac{1}{s-\hat{\lambda}_i}\bfsfc_i\bfsfb_i^T.
\end{equation}
Then $\bfG_r(s)$ is port-Hamiltonian, asymptotically stable, and passive. Moreover, $\bfG_r(s)$ satisfies the necessary condition (\ref{H2optcond1}),  for $\htwo$ optimality: $\bfG_r(s)$ interpolates $\bfG(s)$ at $-\widehat{\lambda}_i$ along the tangent directions $\bfsfb_i$
for $i=1,\ldots,r$.
\end{thm}
{\bf Proof:}
The port-Hamiltonian structure, and consequently passivity, are direct consequences of the construction of $\bfG_r(s)$ as shown in Theorem \ref{thm:origcoor}. The $\htwo$-based tangential interpolation property results from the way the interpolation points $s_i$ and the tangential directions $\bfsfb_i$ are corrected
in Steps 6-c) and 6-d) throughout the iteration. Hence upon convergence; $s_i=-\widehat{\lambda}_i$
and $\bfsfb_i$ is obtained from the residue of  $\bfG_r(s)$ corresponding to $\widehat{\lambda}_i$ as desired.

To prove asymptotic stability we proceed as follows:  Due to Theorem \ref{thm:equivalent}, without loss of generality we can assume that $\bfQ = \bfI$. 
Choose $\bfM$ in \emph{IRKA-PH} so that $\hbfVr$ is orthogonal, i.e. $\hbfVr^T\hbfVr = \bfI_r$
Hence, we obtain $\hbfWr = \hbfVr$. Then, the reduced-order quantities are given by the one-sided projection 
$
\bfA_r = \hbfVr^T \bfA \hbfVr$, and $\bfB_r = \bfC_r^T = \hbfVr^T \bfB.
$
Let $\bfA_r = \bfX_r \mathbf{\Lambda} \bfX_r^{-1}$ be the eigenvalue decomposition of $\bfA_r$ upon convergence of \emph{IRKA-PH} where $\mathbf{\Lambda} = {\rm diag}(\widehat{\lambda}_1,\ldots,\widehat{\lambda}_1)$. Note that  upon convergence,  the tangential directions $\bfsfb_i$ are
obtained by  transposing the rows of the reduced-$\bfB$ matrix of $\bfG_r(s)$  represented in the modal form. 
Define 
$
\bfF_r = [\bfsfb_1,\ldots,\bfsfb_r]^T.
$ 
Hence, upon completion of \emph{IRKA-PH}, we have 
\begin{eqnarray*} 
\bfA_r = \hbfVr^T \bfA \hbfVr = \bfX_r \mathbf{\Lambda} \bfX_r^{-1}~~{\rm and}~~
\bfB_r =  \hbfVr^T \bfB = \bfX_r \bfF_r
\end{eqnarray*}
From the construction of $\bfV_r$ in Step (6)-(e) of \emph{IRKA-PH} and from the fact that 
$s_i = -\widehat{\lambda}_i$ upon convergence, the $i^{\rm th}$ column of $\bfV_r$ satisfies
$(-\widehat{\lambda}_i \bfI - \bfA) \bfv_i = \bfB \bfsfb_i$. Consequently,
 $\bfV_r$ solves a special Sylvester equation
\begin{equation} \label{eqn:sylvesterV}
\bfA \bfV_r + \bfV_r \mathbf{\Lambda}^T + \bfB \bfF_r^T = \mathbf{0}
\end{equation}
where we used the fact $\mathbf{\Lambda} = \mathbf{\Lambda}^T$. Multiply
(\ref{eqn:sylvesterV}) by $\bfX_r^T$ from right and use $\hbfVr = \bfV_r \bfM$ to obtain
\begin{equation} \label{eqn:sylvesterVhat}
\bfA \hbfVr \bfK_r  + \hbfVr \bfK_r \bfA_r^T + \bfB \bfB_r^T = \mathbf{0}
\end{equation}
where $\bfK_r = \bfM^{-1} \bfX_r^T \in \IC^{r \times r}$ is non-singular. Then multiplying (\ref{eqn:sylvesterVhat})
by $\hbfVr^T$ from left yields
\begin{equation} \label{eqn:redlyap}
\bfA_r \bfK_r  + \bfK_r \bfA_r^T + \bfB_r \bfB_r^T = \mathbf{0}
\end{equation}
Since $\bfG_r(s)$ is passive, we know that $\bfA_r$ has no eigenvalues in the open right-half plane.
Hence, to prove  asymptotic stability of $\bfG_r(s)$,
we need to further  show that $\bfA_r$ has no eigenvalues on the imaginary axis either.
Assume to the contrary that  it has, i.e.
$$
\bfz_r^* \bfA_r = \bfz_r^* \lambda,~~{\rm with}~~ \lambda = \jmath \omega.
$$
Note that, $\bfA_r^T \bfz_r = \bar{\lambda}\, \bfz_r = -\lambda\, \bfz_r.$
Multiply (\ref{eqn:redlyap}) by $\bfz_r^*$ from left and by $\bfz_r$ from right
to obtain
\begin{equation} \label{brz}
\bfB_r^T\bfz_r = 0.
\end{equation}
Then, multiplying (\ref{eqn:sylvesterVhat}),  from right, by $\bfz_r$ and using (\ref{brz}) leads to
$$
\bfA \hbfVr \bfK_r \bfz_r = \hbfVr \bfK_r \bfz_r \lambda.
$$
However, this means that $\hbfVr \bfK_r \bfz_r$ is an eigenvector of $\bfA$
with the corresponding eigenvalue $\lambda = \jmath \omega$, contradicting
asymptotically stability of $\bfG(s)$. Therefore, $\bfA_r$ cannot have a
purely imaginary eigenvalue and $\bfG_r(s)$ is asymptotically stable. 
$\Box$\\

The convergence behavior of \emph{IRKA} has been studied in detail in \cite{gugercin2008hmr}; see, also \cite{antoulas2010imr}.  In the vast majority of cases \emph{IRKA} converges rapidly to the desired first-order conditions independent of initialization. Effective initialization strategies have been proposed in  \cite{gugercin2008hmr} as well, although random initialization often performs very well. 
We illustrate the robustness of \emph{IRKA-PH} with regard to initialization in Section \ref{sec:examples}.
\begin{rem} Let $\PH(r)$  denote the set of port-Hamiltonian systems with state-space dimension $r$ 
as in \eqref{eq:lphs_rom} and consider the problem: 
\begin{equation} \label{h2optPH}
 \left\| \bfG - \bfG_r \right\|_{\htwo} =  \min_{ \mbox{\tiny{$
\widetilde{\bfG}_{r}\!\in\!\PH(r)$}} } \left\| \bfG - \widetilde{\bfG}_{r} \right\|_{\htwo}
\end{equation}
where $\bfG(s)$ and $\bfG_r(s)$  are both port-Hamiltonian.  This is the problem that we would prefer to solve and it is a topic of current research.
Note that although \emph{IRKA-PH} generates a reduced-order port-Hamiltonian model, $\bfG_r(s)$, satisfying a first-order necessary condition for $\htwo$ optimality,  this condition is not a necessary condition for $\bfG_r(s)$ to solve (\ref{h2optPH}), that is, for $\bfG_r(s)$ to be the \emph{optimal} reduced-order port-Hamiltonian system approximation to the port-Hamiltonian system $\bfG(s)$.  Nonetheless, we find that reduced-order models produced by \emph{IRKA-PH} have much superior $\htwo$ performance compared to other approaches. This is illustrated in Section \ref{sec:examples}. 
\end{rem}

\begin{rem} \label{PHandH2Opt}
\emph{IRKA-PH}, as described in Algorithm \ref{phirka}, produces a reduced order port-Hamiltonian model $\bfG_r(s)$, which satisfies the first-order necessary condition of $\htwo$-optimality (\ref{H2optcond1}).  The remaining two conditions for $\htwo$-optimality, (\ref{H2optcond2}) and (\ref{H2optcond3}), will be satisfied if in addition:
\begin{eqnarray}
\label{eq:algebraic_condition}
\begin{array}{l}
\rm{Range}[(\widehat{\lambda}_1\bfI+(\bfJ - \bfR)\bfQ)^{-1}\bfB\bfsfb_1,\dots,\\
\hspace{30mm}(\widehat{\lambda}_r\bfI+(\bfJ - \bfR)\bfQ)^{-1}\bfB\bfsfb_r] =\\
\rm{Range}[(\widehat{\lambda}_1\bfI+(\bfJ - \bfR)^T\bfQ)^{-1}\bfB\bfsfc_1,\dots,\\
\hspace{30mm}(\widehat{\lambda}_r\bfI+(\bfJ - \bfR)^T\bfQ)^{-1}\bfB\bfsfc_r],
\end{array}
\end{eqnarray}
where $\widehat{\lambda}_i$ are reduced system poles and $\bfsfb_i,\bfsfc_i$ are the corresponding (vector) residues in the expansion (\ref{partFracExp}) for $i=1,\ldots,r$.   The condition (\ref{eq:algebraic_condition}) allows retention of port-Hamiltonian structure while enforcing the $\Htwo$-optimal bitangential interpolation conditions (\ref{H2optcond1})-(\ref{H2optcond3}).     Typically (\ref{eq:algebraic_condition}) will not hold, however the special case when 
$\bfJ=\mathbf{0}$ will lead to (\ref{eq:algebraic_condition}) with $\bfsfb_i=\bfsfc_i$ implying satisfaction of $\Htwo$ optimality conditions for the reduced system. 

\end{rem}

\section{Numerical Examples}  \label{sec:examples}

In this section, we illustrate the theoretical discussion on three port-Hamiltonian systems. The first two systems are of modest dimension so that we can compute all the norms (such as $\htwo$ and $\Hinf$ error norms) explicitly for several $r$ values for many different methods to compare them to the fullest extent. Then, to illustrate that the proposed method can be easily applied in the large-scale settings as intended, we use a large-scale model as the third example.

\subsection{MIMO Mass-Spring-Damper system}

\begin{figure}[t]
\centering
\includegraphics[width=8cm]{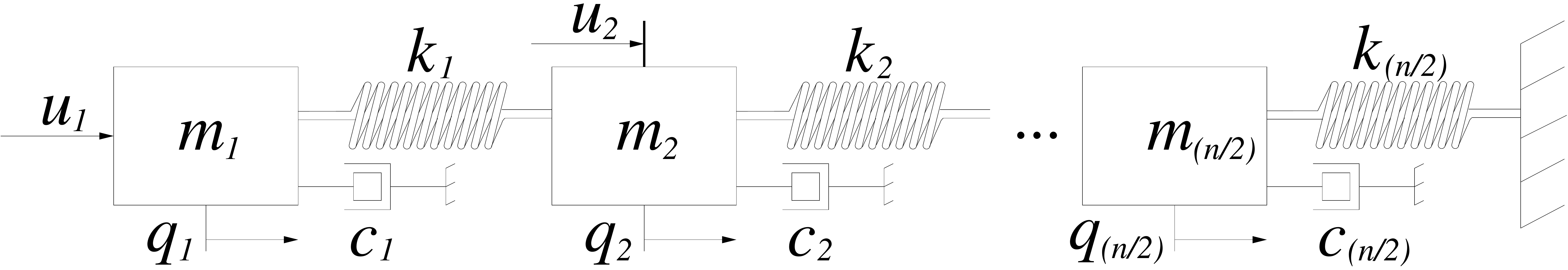}
\caption{Mass-Spring-Damper system}
\label{fig:msd}
\end{figure}

The full-order model is the the Mass-Spring-Damper system shown in Fig. \ref{fig:msd} with masses $m_i$, spring constants $k_i$ and damping constants $c_i \geq 0, \ i = 1, \dots, n/2$. $q_i$ is the displacement of the mass $m_i$. The inputs $u_1, u_2$ are the external forces applied to the first two masses $m_1, m_2$. The port-Hamiltonian outputs $y_1, y_2$ are the velocities of masses $m_1, m_2$. The state variables are as follows: $x_1$ is the displacement $q_1$ of the first mass $m_1$, $x_2$ is the momentum $p_1$ of the first mass $m_1$, $x_3$ is the displacement $q_2$ of the second mass $m_2$, $x_4$ is the momentum $p_2$ of the second mass $m_2$, etc.


A minimal realization of this port-Hamiltonian system for order $n = 6$ corresponding to three masses, three springs and three dampers is
\begin{eqnarray*}
\bfB= \begin{bmatrix}~0~1~0~0~0~0~ \\~0~0~0~1~0~0~\end{bmatrix}^T  , \
\bfC =\begin{bmatrix}~0~\frac{1}{m_1}~0~0~0~0~\\~0~0~0~\frac{1}{m_2}~0~0~\end{bmatrix},
\end{eqnarray*}
\vspace*{-0.5cm}
\begin{eqnarray*}
\bfJ=
\left[ \begin{array}{cccccc}
0 &  1 &  0 &  0 &  0 & 0 \\
-1 &  0 &  0 & 0 &  0 & 0 \\
0 &  0 &  0 & 1 &  0 & 0 \\
0 &  0 &  -1 & 0 &  0 & 0 \\
0 &  0 &  0 & 0 &  0 & 1 \\
0 &  0 &  0  &  0 &  -1 & 0
\end{array}  \right] ,~~
\bfR=
\left[ \begin{array}{cccccc}
0 &  0 &  0 &  0 &  0 & 0 \\
0      &  c_1 &  0 & 0 &  0 & 0 \\
0 &  0 &  0 & 0 &  0 & 0 \\
0 &  0 &  0 & c_2 &  0 & 0 \\
0 &  0 &  0 & 0 &  0 & 0 \\
0      &  0 &  0      &  0 &  0 & c_3
\end{array}  \right] , 
\end{eqnarray*}
{\rm and}
\begin{eqnarray*}
\bfQ = \left[\begin{array}{cccccc}
k_1 & 0 & -k_1 & 0  & 0  & 0   \\
0 & \frac{1}{m_1} & 0 & 0  & 0  & 0 \\
-k_1 & 0 & k_1 +k_2 & 0 & -k_2  & 0   \\
0 & 0 & 0 & \frac{1}{m_2}  & 0  & 0  \\
0 & 0 & -k_2 & 0  & k_2+k_3  & 0 \\
0 & 0 & 0 & 0  & 0  & \frac{1}{m_3}
\end{array}\right].
\end{eqnarray*}

Then the $\bfA$ matrix can be obtained  as
\begin{eqnarray*}
\begin{array}{l}
\bfA = (\bfJ - \bfR)\bfQ = \\
\left[\begin{array}{cccccc}
0  & \frac{1}{m_1} & 0 & 0  & 0  & 0   \\
-k_1 & -\frac{c_1}{m_1} & k_1 & 0  & 0  & 0    \\
0 & 0 & 0 & \frac{1}{m_2}  & 0  & 0  \\
k_1 & 0 & -k_1 -k_2 & -\frac{c_2}{m_2}  & k_2  & 0 \\
0 & 0 & 0 & 0  & 0  & \frac{1}{m_3}  \\
0 & 0 & k_2 & 0  & -k_2-k_3  & -\frac{c_3}{m_3}
\end{array}\right].
\end{array}
\end{eqnarray*}


Adding another mass with a spring and a damper would increase the dimension of the system by two. This will lead to
a zero entry in the $(n -1, n -1)$ position and an entry of $-c_{n/2}/m_{n/2}$ in the $(n, n)$ position. The superdiagonal of $\bfA$ will have $k_{n/2 - 1}$ in the $(n -2, n -1)$ position and $1/m_{n/2}$ in the $(n -1, n)$ position. The subdiagonal of $\bfA$ will have $0$ in the $(n - 1, n - 2)$ position and $-k_{n/2 - 1} - k_{n/2}$ in the $(n, n -1)$ position. Additionally $\bfA$ will have $k_{n/2 - 1}$ in the $(n, n -3)$ position.


We used a 100-dimensional Mass-Spring-Damper system with  $m_i = 4 ,\ k_i =  4$, and $ c_i = 1$.
We compare three different methods: (1) The  proposed method in Algorithm \ref{phirka} denoted by \emph{IRKA-PH}
(2) The effort-constraint method \eqref{eq:RedSysEfConstr} of 
(\cite{Polyuga_vdSchaft2011_Effort_flow,vdSchaft_Polyuga_09cdc,Polyuga_2010_thesis,PolvdSc08,PolvdSc_2011_lecture_notes})
denoted by \emph{Eff. Bal.}, and (3) One step interpolatory model reduction without the $\htwo$ iteration, denoted by
\emph{1-step-Intrpl}. The reason for including \emph{1-step-Intrpl} is as follows: We choose a set of interpolation points and tangential directions; then using the interpolatory projection in Theorem \ref{thm:origcoor}, we obtain a reduced model. Then, we use the same interpolation points and directions as an initialization technique for \emph{IRKA-PH}. Then, we are able to illustrate starting from the same interpolation data, how
 \emph{IRKA-PH} corrects the interpolation points and directions, and how much it improves the $\htwo$ and $\Hinf$ behavior through out this iterative process.

\begin{figure}[t]
\centering
\includegraphics[width=8cm]{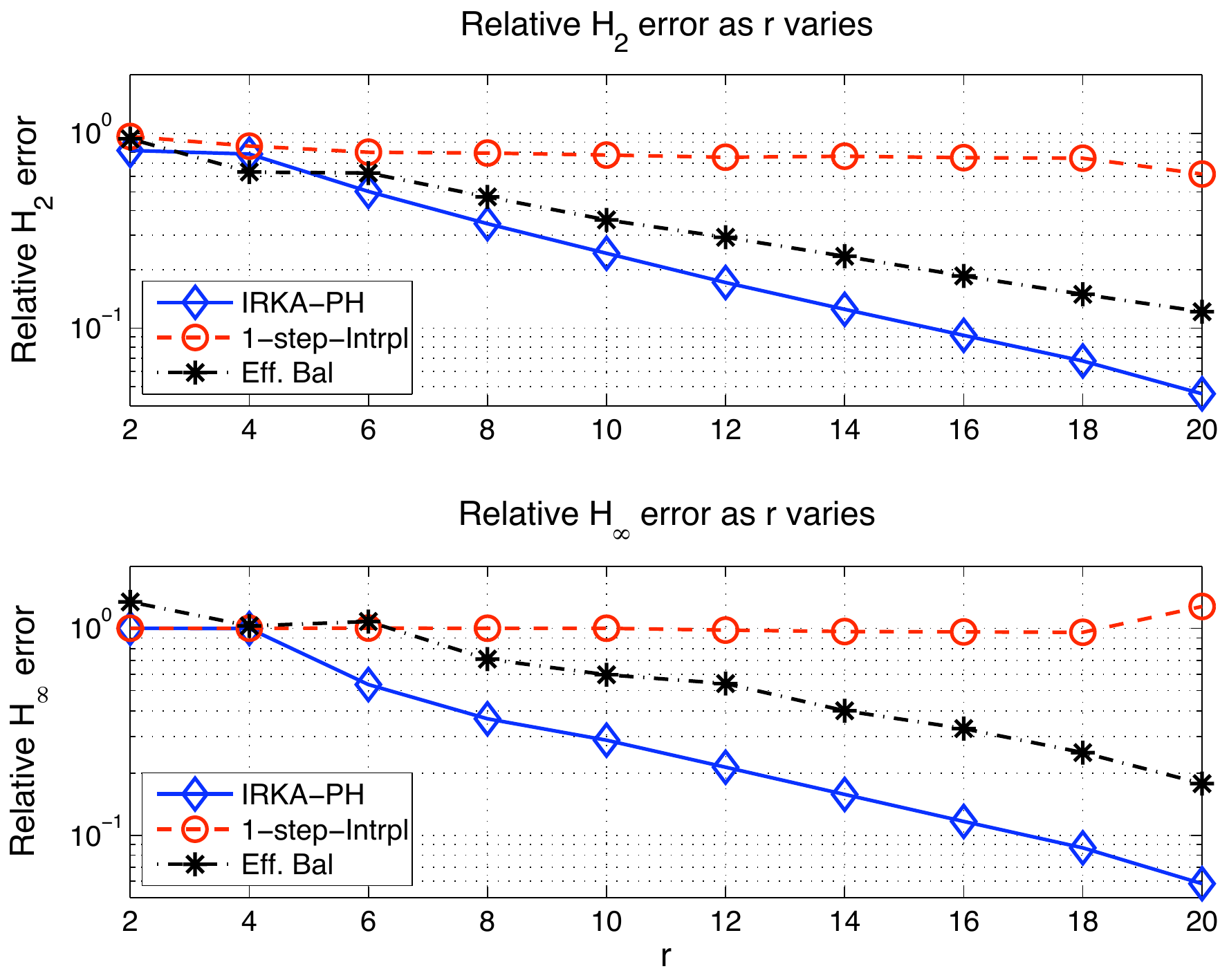}
\caption{Relative $\htwo$ and $\hinf$ norms for Mass-Spring-Damper System}
\label{fig:msd_h2_hinf}
\end{figure}

\begin{figure}[t]
\centering
\includegraphics[width=8cm]{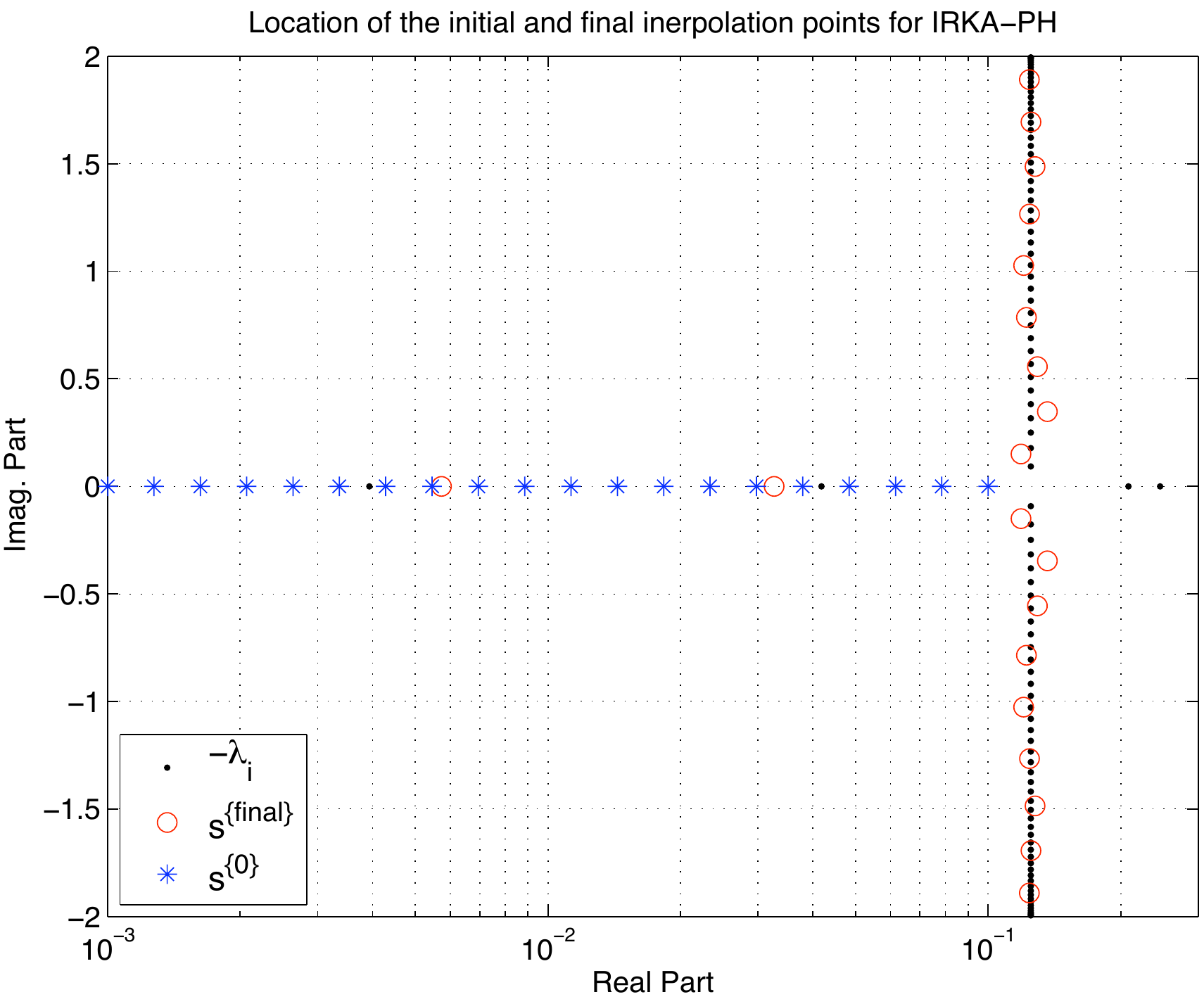}
\caption{Initial and converged interpolation points}
\label{fig:msd_shifts}
\end{figure}

Using all three methods, we reduce the order from $r=2$ to $r=20$ with increments of two. For \emph{IRKA-PH},
initial interpolation points are chosen as logarithmically spaced points between $10^{-3}$ and $10^{-1}$; and the corresponding directions are the dominant  right singular vectors of the transfer function at each interpolation point (a small $2 \times 2$
singular value problem). These same points and directions are also used for \emph{1-step-Intrpl}.
The resulting relative $\htwo$ and $\hinf$ error norms
for each order $r$ are illustrated in Figure \ref{fig:msd_h2_hinf}. Several important conclusions are in order: First of all,
with respect to the both $\htwo$ and $\hinf$ norm, \emph{IRKA-PH} significantly outperforms the other methods. As the figure illustrates, for \emph{1-step-Intrpl}, the performance hardly improves as $r$ increases unlike \emph{IRKA-PH} for which both the $\htwo$ and $\hinf$ errors decay consistently. This shows the superiority of \emph{IRKA-PH} over
 \emph{1-step-Intrpl}. The initial interpolation point and direction selection does not yield a satisfactory interpolatory reduced-order model; however instead of searching some other interpolation data in an {\it ad hoc} way, the proposed method automatically corrects the interpolation points through out the iteration and yields a significantly smaller error norm. To see this effect more clearly, we plot the initial interpolation point selection, denoted by ${\rm s^{\{0\}}}$, and the final/converged interpolation points, denoted by ${\rm s^{\{final\}}}$, in Figure \ref{fig:msd_shifts} together with the
 mirror images of the original system poles, denoted by $-\lambda_i$. As this figure illustrates, even starting with this logarithmically placed points on the real line, \emph{IRKA-PH} iteratively corrects the points so that upon convergence
they automatically align themselves around the mirror images of the original systems poles.
They don't necessarily exactly match the mirror images of the original poles but rather they align themselves in a way that overall they have a good representation of the original spectrum. This is similar to the eigenvalue computations where the Ritz values  provide better information about the overall spectrum than just a subset of the exact eigenvalues.  We further note that the full-order poles are computed only to obtain this figure and are not needed by  \emph{IRKA-PH} .

The second observation from Figure \ref{fig:msd_h2_hinf} is that in comparison to \emph{Eff. Bal.}, \emph{IRKA-PH} achieves the better performance with less computational  effort; the main cost is sparse linear solves and no dense matrix operations are needed unlike the balancing-based approaches where Lyapunov equations need to be solved.
 We also note that even though the proposed method is $\htwo$-based, it produces a satisfactory $\hinf$ performance as well. This is not surprising as \cite{gugercin2008hmr} discusses \emph{IRKA} usually leads to high fidelity  $\hinf$ performance as well. 


\begin{figure}[t]
\centering
\includegraphics[width=8cm]{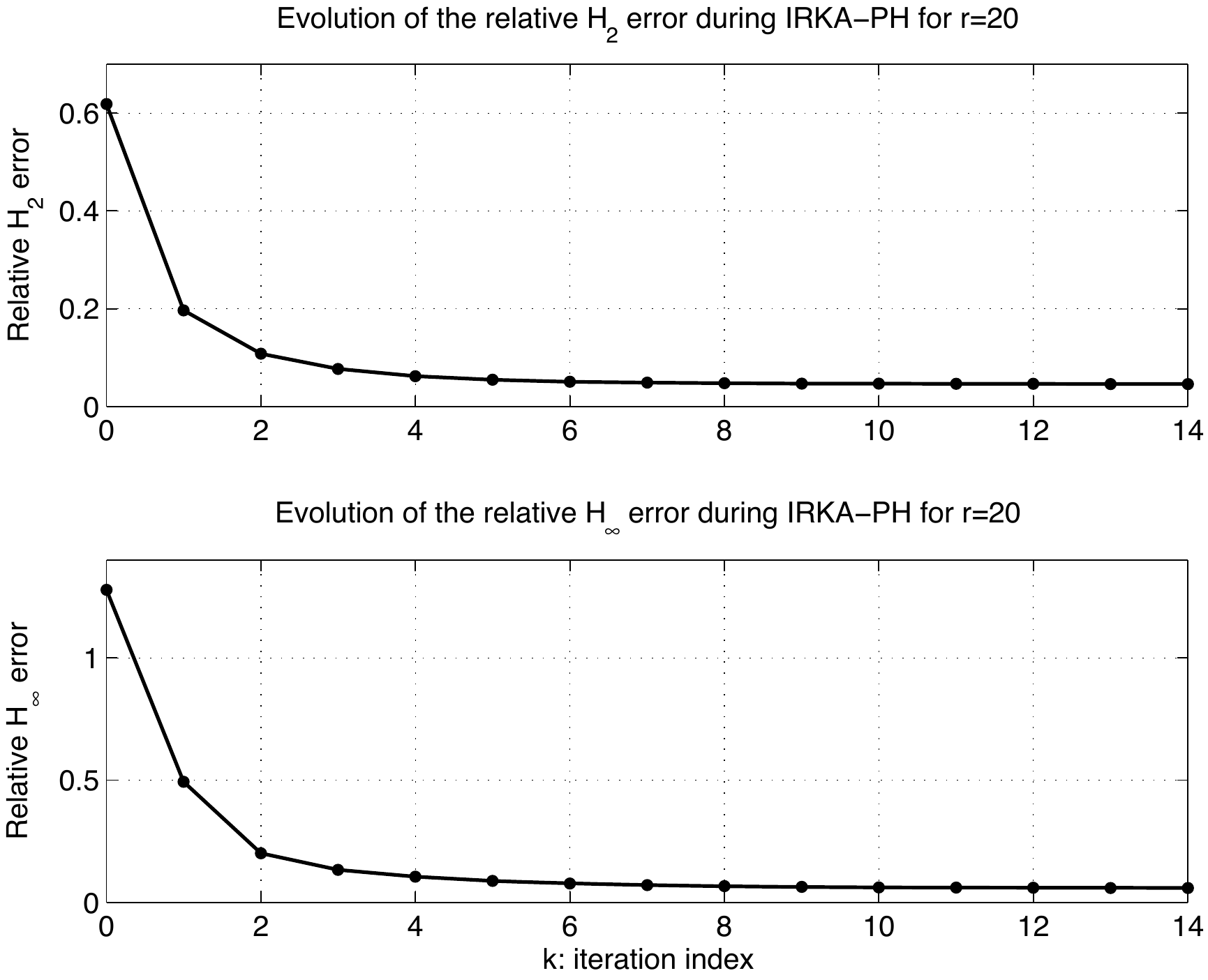}
\caption{Relative $\htwo$ and $\hinf$ norms for Mass-Spring-Damper System}
\label{fig:msd_h2hinf_evol}
\end{figure}

Before we conclude this example, we further investigate the $r=20$ case. We show, in Fig \ref{fig:msd_h2hinf_evol}, how the $\htwo$ and $\hinf$ errors evolve
during \emph{IRKA-PH} for $r=20$, the largest reduced-order. The figure reveals convergence within seven or eight steps.  The large initial relative errors are reduced drastically already after two steps of the iteration.  This illustrates the effectiveness of the proposed method.

Next, we investigate the effect of different initializations on the performance of  \emph{IRKA-PH}. For brevity, we only illustrate the $r=20$ case.   
We make $5$ different initializations and denote by $\mathcal{S}_0^{\{j\}}$ the set of initial interpolation points corresponding to the $j^{\rm th}$ selection. $\mathcal{S}_0^{\{1\}}$ will be the same as what was used earlier, i.e., $20$ points logarithmically spaced between $10^{-3}$ and $10^{-1}$. For $\mathcal{S}_0^{\{2\}}$, we choose $20$ points logarithmically spaced between $-10^{-5}$ and $-10^{-2}$.  Note that this is a poor selection since these initial interpolation points lie in the left-half plane in proximity of the poles. We consciously make this (bad) choice to see the effect on convergence.  
For $\mathcal{S}_0^{\{3\}}$, we choose  complex points in the right half-plane with real parts that are logarithmically spaced between $10^{-6}$ and $10^0$' and imaginary parts that are logarithmically spaced between $10^{-3}$ and $10^{-1}$.  The set is chosen to be closed under conjugation.  These points are arbitrarily selected and have nothing to do with the spectrum of $\bfA$. For $\mathcal{S}_0^{\{4\}}$, we make the situation even worse than for $\mathcal{S}_0^{\{2\}}$. We choose $20$ poles of the original system, $\bfG(s)$, and perturb them by $0.1\%$ to obtain our starting points. This is a very bad selection in two respects: 1) The interpolation points lie in the left-half plane;  2) They are extremely close to system poles and make the linear system $(s_i \bfE - \bfA) \bfv_i = \bfB \bfsfb_i$ very poorly conditioned.  Finally, for $\mathcal{S}_0^{\{5\}}$ we choose, once again, $r=20$ original systems poles, but this time reflect them across the imaginary axis to obtain initialization points. Once interpolation points are chosen,  associated directions are taken to be the dominant  right singular vectors of the transfer function at each interpolation point as we did previously. 

\begin{figure}[t]
\centering
\includegraphics[width=8cm]{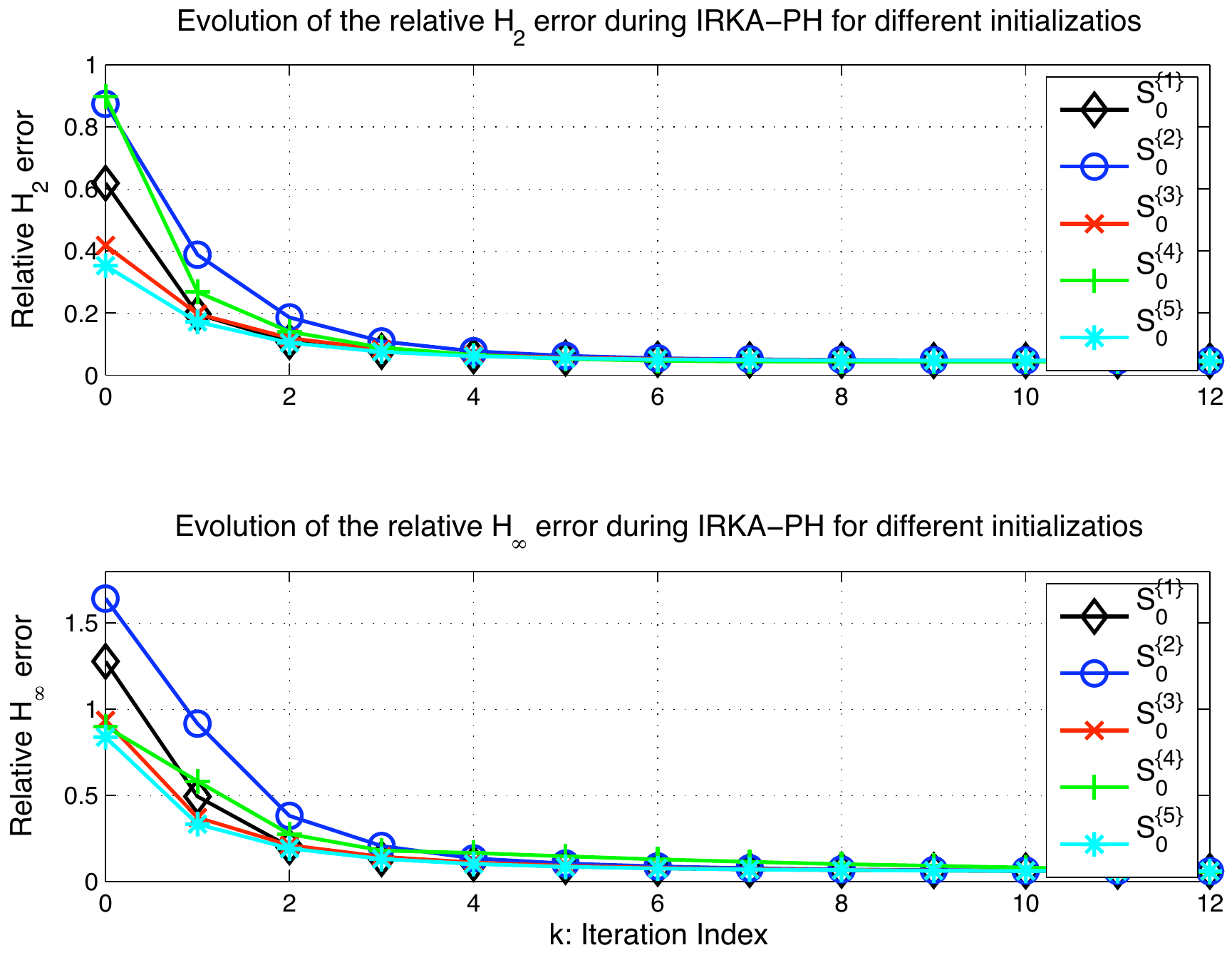}
\caption{Relative $\htwo$ and $\hinf$ norms for Different Initializations}
\label{fig:msd_h2hinf_evol_diff_init}
\end{figure}

Figure \ref{fig:msd_h2hinf_evol_diff_init} shows the evolution of 
the relative $\htwo$ and $\hinf$ errors during \emph{IRKA-PH} for these $5$ different selections. 
In all cases,  \emph{IRKA-PH} converges to the same reduced-model in almost the same number of steps. As expected,  $\mathcal{S}_0^{\{2\}}$ and $\mathcal{S}_0^{\{4\}}$ are the worst initializations, starting with relative errors bigger than $1$ in the $\hinf$ norm. However, the algorithm successfully corrects these points and drives them towards high-fidelity interpolation points and tangent directions. Even though $\mathcal{S}_0^{\{5\}}$  seems to be the best initialization --- it starts with the lowest initial error --- it converges to the same reduced-model as the iteration that started with $\mathcal{S}_0^{\{1\}}$, and also in the same number of steps.  $\mathcal{S}_0^{\{1\}}$ achieves this without the need for original poles. This numerical evidence is in support of the robustness of \emph{IRKA-PH}. The algorithm successfully corrects the bad initializations. Indeed, \emph{IRKA-PH}
is expected to be more robust and to converge faster than even the original \emph{IRKA}, which already depicts fast convergence behavior. The reasons are twofold. Here, unlike in \emph{IRKA}, regardless of the initialization, every intermediate reduced-model is stable; hence the interpolations points never
 lie in the left-half plane. This already smoothens the convergence behavior. Moreover, unlike \emph{IRKA} which uses an oblique projector, \emph{IRKA-PH} is theoretically equivalent to using an orthogonal projector (with respect to an inner product weighted by $\bfQ$)) This in turn supports the robustness of the algorithm. 
 
For all these reasons, we expect \emph{IRKA-PH} to be a numerically robust and rapidly converging iteration.
Throughout our numerical experiments using a wide variety of different initializations, we have \emph{never} experienced a convergence failure of  \emph{IRKA-PH}.  \emph{IRKA-PH}  always converged to the same reduced-model regardless of initialization; even a search for a counterexample spanning thousands of trials failed to produce even a single case of either convergence failure or convergence to a different model in the $r=20$ case.  Indeed this was true throughout the range $8 \leq r \leq 20$;  
\emph{IRKA-PH} converged to the same reduced-model for many different initializations for $r=8:2:20$. Only for 
$r=2$, $r=4$, and $r=6$ and then only after several trials, were we able to make  \emph{IRKA-PH} converge to a different reduced-model  (only \emph{one} different model) than what we had 
before. The performance of this new model was slightly better than what we had before; hence it would further  support the success of \emph{IRKA-PH}. However, as our simple initialization resulted in the same accuracy for every $r=8:2:20$ and did only slightly worse for 
$r=2,4,6$, there is no need to update our earlier results. These numerical experiments once 
more strongly points toward the robustness of \emph{IRKA-PH}.

\begin{figure}[t]
\centering
\epsfxsize= 3.4in
\epsfysize = 4in
\epsffile{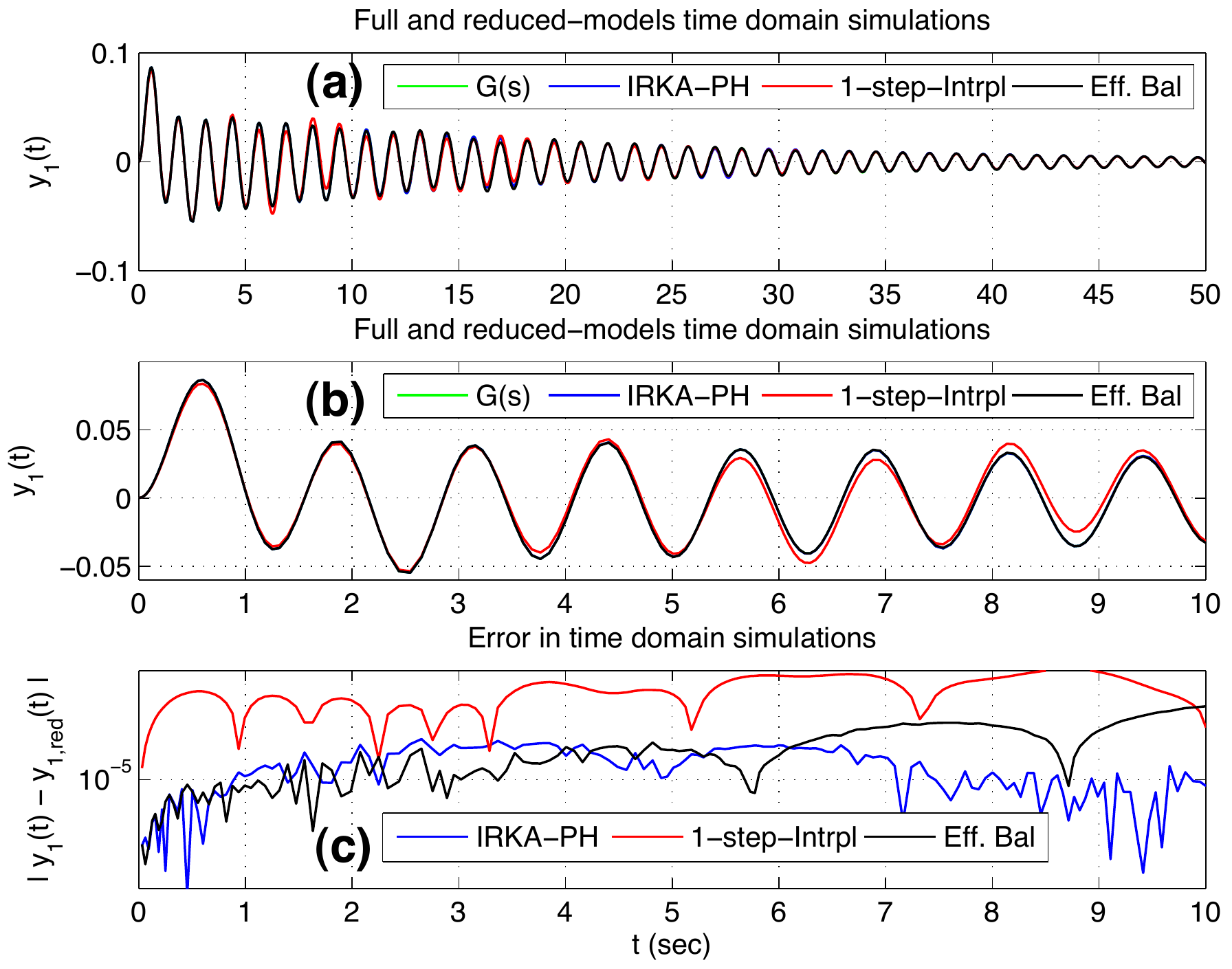}
\caption{Time domain simulations for $r=20$}
\label{fig:msd_r_20_time}
\end{figure}

We present further comparison for the $r=20$ case, by comparing time domain simulations for the full and reduced-models. Recall that $\bfG(s)$ has $2$ inputs and $2$ outputs. To make the time-domain illustrations simpler, we only compare the outputs of the subsystem relating the first-input ($\bfu_1$) to the first-output ($\bfy_1$). The results for the other $3$ subsystems depict the same behavior.
As the input, we choose a decaying sinusoid 
$\bfu_1(t) = e^{-0.05 t} \sin{(5t)}$. We run the simulations for $T=50$ seconds as it takes long time for the response to decay due to the poles  close to the imaginary axis. The results are depicted in Figure \ref{fig:msd_r_20_time}. In Figure \ref{fig:msd_r_20_time}-(a), we plot the simulation results for the whole time interval. To give a better illustration, in Figure \ref{fig:msd_r_20_time}-(b), we zoom into the time interval $[0,10]$ seconds which illustrates that 
1-Step-Intrpl leads to the largest deviation. 
In Figure \ref{fig:msd_r_20_time}-(c), we give the absolute value of the error between the true system output and the reduced ones. As it is clear from this figure and as expected from 
the earlier analysis as shown in Figure \ref{fig:msd_h2_hinf}, \emph{IRKA} yields the smallest deviation. 
The maximum absolute values of the output errors, i.e. $\max | \bfy_1(t) - \bfy_{1,red}(t) |$ where $\bfy_{1,red}(t)$ denotes the first-output  due to the reduced-models
 are 
\begin{center}
\begin{tabular}{c|c|c}
  IRKA-PH & 1-Step-Intrpl  & Effort-Bal  \\ \hline
  $1.31 \times 10^{-3}$ &   $1.09 \times 10^{-2}$ & $3.96 \times 10^{-3}$ \\
\end{tabular}
\end{center}

\subsection{MIMO port-Hamiltonian Ladder Network}
\begin{figure}[hhh]
\centering
\includegraphics[width=8.5cm]{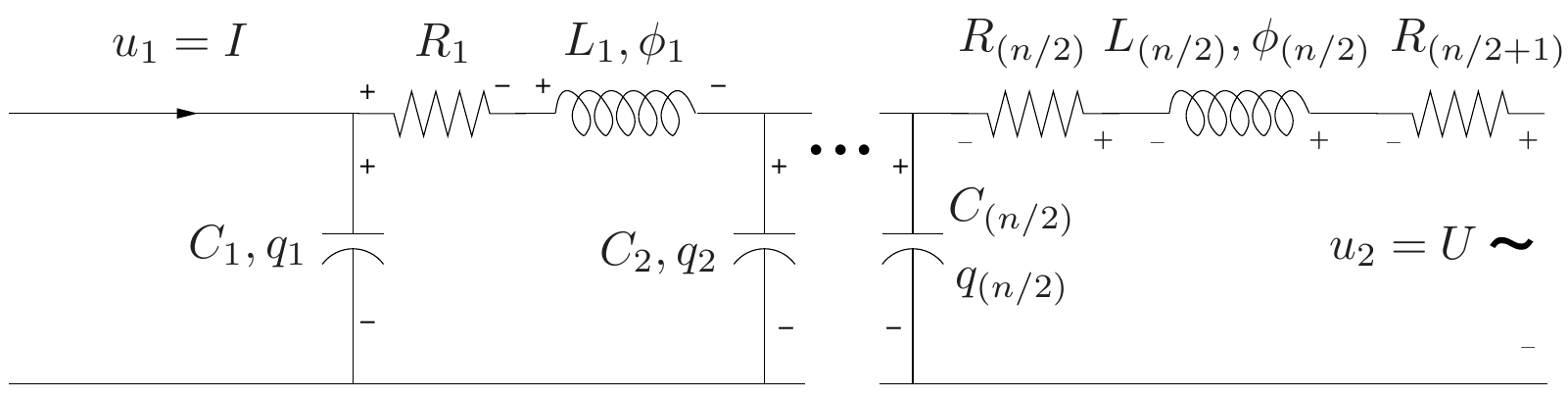}
\caption{MIMO Ladder Network}
\label{fig:mimo1}
\end{figure}

As a second MIMO  port-Hamiltonian system, we consider an $n$-dimensional
Ladder Network as shown in Fig. \ref{fig:mimo1}.

We take the current $I$ on the left side and
the voltage $U$ on the right side of the Ladder Network as the inputs. The
port-Hamiltonian outputs are the voltage over the first capacitor $U_{C_1}$ and the
current through the last inductor $I_{L_{n/2}}$. The state variables are as follows:
$x_1$ is the charge $q_1$ of $C_1$, $x_2$ is the flux $\phi_1$ of $L_1$, $x_3$ is
the charge $q_2$ of $C_2$, $x_4$ is the flux $\phi_2$ of $L_2$, etc. The directions
chosen for the internal currents of the Network are shown by plus- and minus-signs
in Fig  \ref{fig:mimo1}. A minimal realization of this port-Hamiltonian Ladder Network
for order $n = 4$ is given in Example \ref{exm:ladder_x}.

Adding another $LC$ pair to the network, which would correspond to an increase of
the dimension of the model by two, will modify the system matrices as follows:
The subdiagonal of the matrix $\bfA$ will contain additionally
$L_{n/2-1}^{-1}, -C_{n/2}^{-1}$ with the plus-sign in the $(n/2+1, n/2)$ position. The
superdiagonal of $\bfA$ will contain $-C_{n/2}^{-1},L_{n/2}^{-1}$ with the
minus-sign in the $(n/2, n/2+1)$ position. Furthermore, the main diagonal of $\bfA$
will have $-\frac{R_{n/2-1}}{L_{n/2-1}}$ in the $(n-2,n-2)$ position, zero in the $(n-1,n-1)$ position, and
$-\frac{R_{n/2}+R_{n/2+1}}{L_{n/2}}$ in the $(n,n)$ position. $\bfB$ and $\bfC$ matrices will be of
the similar structure as in \eqref{eq:ladder_ABC} and \eqref{eq:ladder_ABCC}. The matrix $\bfB$ will have ones in the
$(1,1)$ and $(n,2)$ positions and zeros in the other positions. The only nonzero
entries of $\bfC$ are $\frac{1}{C_1}$ in the $(1,1)$ position and $\frac{1}{L_{n/2}}$ in
the $ (2,n)$ position.

Here, we consider  the 100-dimensional full order port-Hamiltonian network with  $C_i=0.1$ and $L_i=0.1$ and $R_i= 3$ for $i=1,\ldots,50$ and $R_{51}=1$. For this model, we compare \emph{IRKA-PH} with both 
\emph{Eff. Bal.} and
{\it regular balanced truncation} (denoted by \emph{Reg. Bal.}).  We would like to emphasize that 
\emph{Reg. Bal.}
 does not preserve the port-Hamiltonian structure.
We choose to include   \emph{Reg. Bal.} in our comparison to better
 illustrate the effectiveness of our proposed method by showing that  \emph{IRKA-PH} can
 perform as good as or sometimes 
  even better than \emph{Reg. Bal.} which is known to yield high-fidelity $\Hinf$ and $\htwo$ performance (though it is not constrained to preserve port-Hamiltonian structure).

\begin{figure}[t]
\centering
\includegraphics[width=8cm]{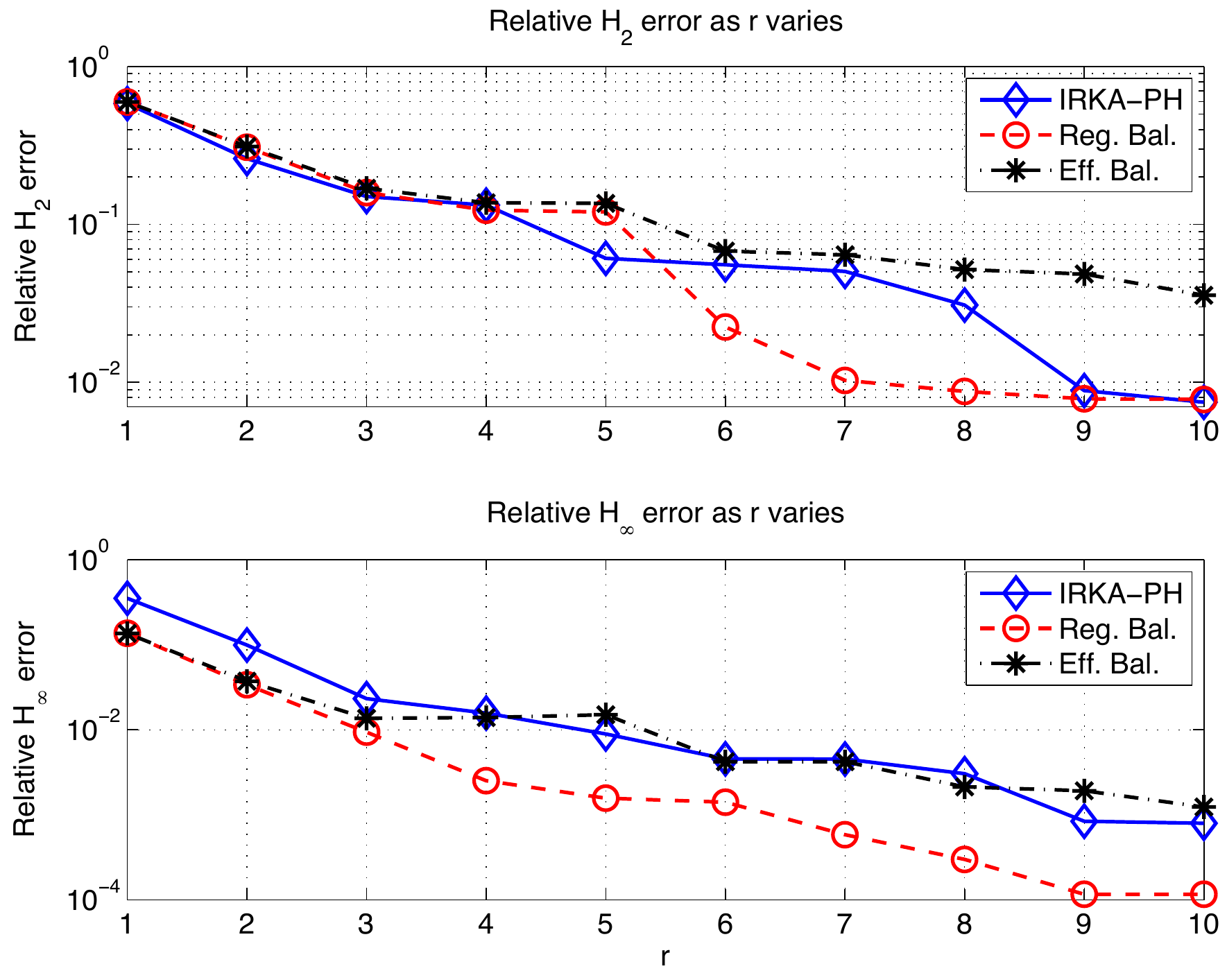}
\caption{Evolution of the relative $\htwo$ and $\hinf$ norms}
\label{fig:ladder_damped_h2_hinf}
\end{figure}

We reduce the order from $r=1$ to $r=10$ in increments of one.
The resulting relative $\htwo$ and $\hinf$ error norms are depicted in
 Figure \ref{fig:ladder_damped_h2_hinf}.
 The $\htwo$-based nature of our proposed method
 is clear. \emph{IRKA-PH} outperforms \emph{Eff. Bal.} for every $r$. Indeed, what is interesting to observe is that \emph{IRKA-PH} is better than even  \emph{Reg. Bal.} for each $r=1,\ldots,5$. Even though 
 for $r=6,7,8$  \emph{Reg. Bal.} is better, for the last two $r$ values, \emph{IRKA-PH} is able to match 
   \emph{Reg. Bal}. Note that our proposed method achieves this performance while still preserving structure without any need for solving Lyapunov equations; the only cost is sparse linear solves. 
  As expected, \emph{Reg. Bal.} is the best in terms of $\Hinf$ performance. Indeed,  it is tailored towards $\Hinf$ error reduction and is not constrained to preserve structure. However, the $\Htwo$-based 
 \emph{IRKA-PH}  performs as well as \emph{Eff. Bal.} in terms of $\Hinf$ error norm. This once more shows that, similar to \emph{IRKA},   \emph{IRKA-PH} provides high-fidelity 
  not only in the $\Htwo$ norm but also  in the $\Hinf$ norm. 
 
\begin{figure}[t]
\centering
\includegraphics[width=8cm]{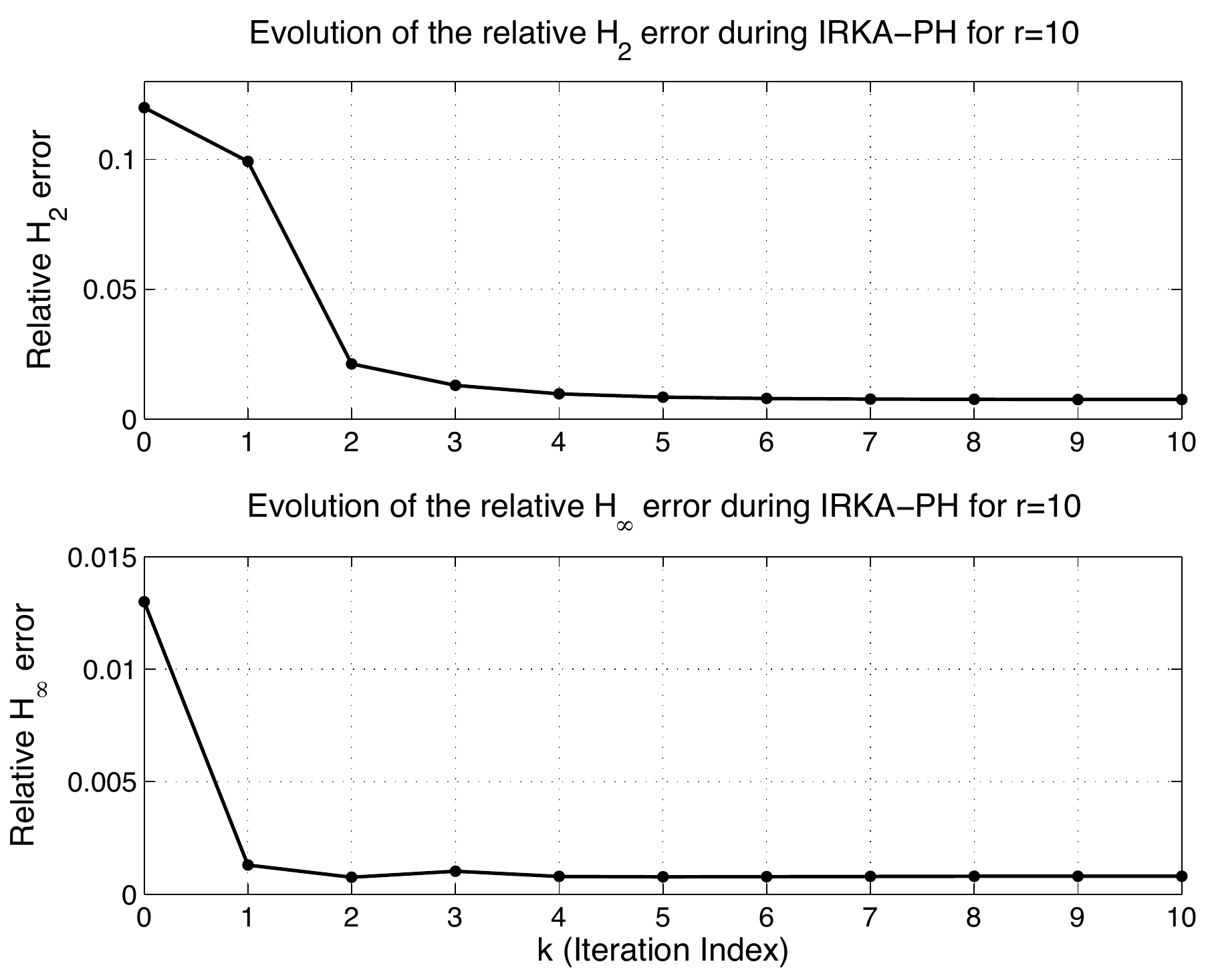}
\caption{Evolution of the relative $\htwo$ and $\hinf$ norms}
\label{fig:ladder_damped_h2hinf_evol}
\end{figure}

Similar to the previous example, we illustrate the convergence behavior of \emph{IRKA-PH}, in this case for $r=10$, in Figure \ref{fig:ladder_damped_h2hinf_evol}. In this case, the convergence is even faster than the previous example with  the algorithm converging after five to six steps. For this model,
we have initialized the interpolation points for \emph{IRKA-PH}  arbitrarily as logarithmically spaced points between $10^{-2}$ and $10^1$. As before, the initial tangent directions are chosen as the leading  right singular vector of $\bfG(s)$ evaluated at the chosen interpolation points.

\begin{figure}[t]
\centering
\includegraphics[width=8cm]{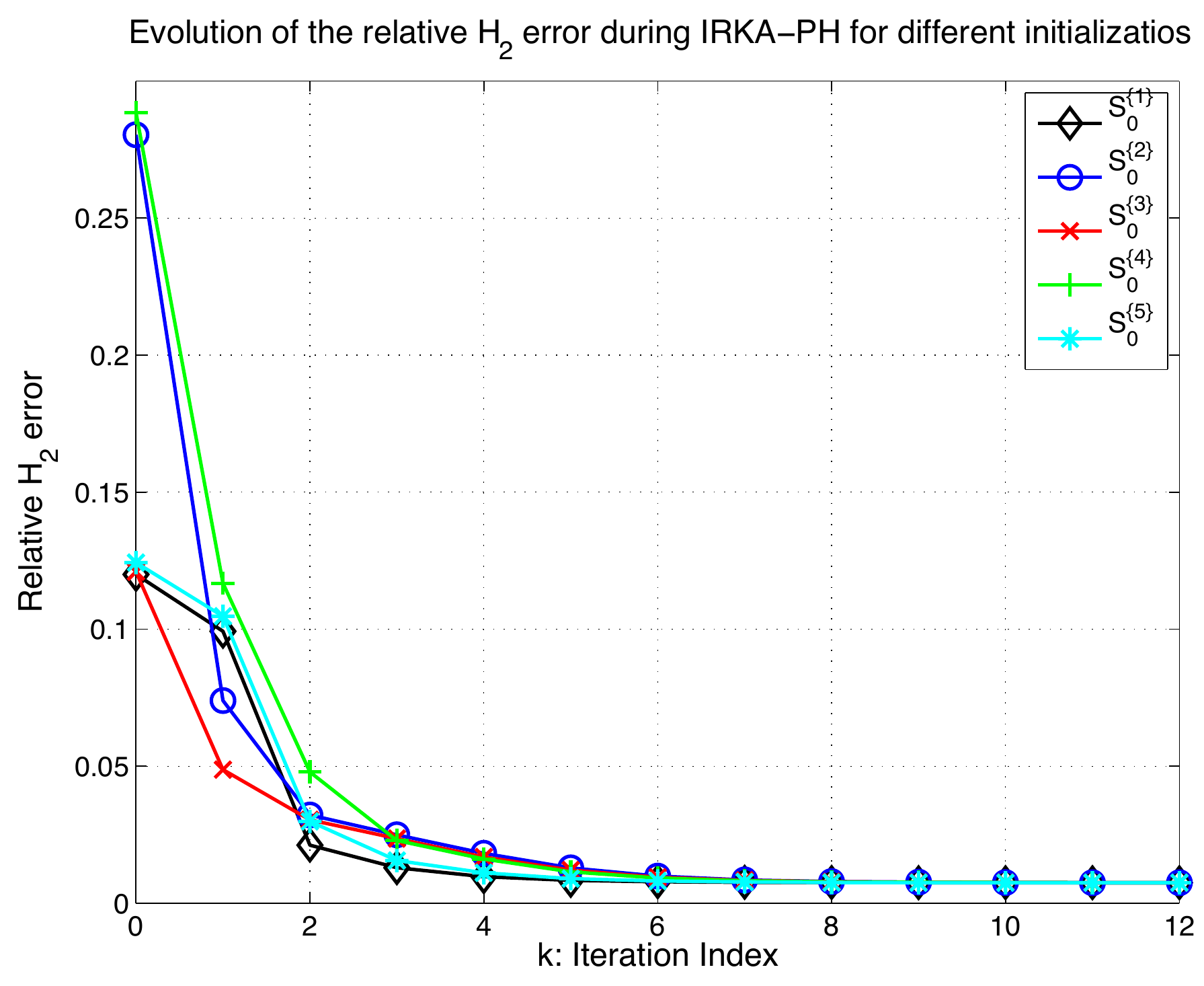}
\caption{Relative $\htwo$  norms for Different Initializations}
\label{fig:ln_damped_h2_evol_diff_init}
\end{figure}

As we did in the previous example, we experiment with different initialization strategies for \emph{IRKA-PH}. For  brevity, we only illustrate the $r=10$ case.  We consider  $5$ different initializations, $\mathcal{S}_0^{\{j\}}$ for $j=1,\ldots,5$ for the interpolation points. $\mathcal{S}_0^{\{1\}}$ is what used earlier. For $\mathcal{S}_0^{\{2\}}$,
we choose $r=10$ original systems poles and perturb them by $0.1\%$  as the starting points.
This is a very poor choice since the interpolation points are very close to the system poles.
For $\mathcal{S}_0^{\{3\}}$, we choose $10$ points logarithmically spaced between $-10^{-5}$ and $-10^{-1}$. Once again, these interpolation points are also in the left-half plane; a strategy one would usually avoid. $\mathcal{S}_0^{\{4\}}$ correspond to choosing $10$ points logarithmically spaced between $10^{-8}$ and $10^{-4}$. 
And finally,  for $\mathcal{S}_0^{\{5\}}$, we choose some arbitrary complex numbers where
the real parts are  distributed between $10^{-5}$ and $10^{-2}$ and
 the imaginary parts are  distributed between $10^{-2}$ and $10^{0}$ together with their complex conjugate pairs. 
Once the starting points are chosen, the corresponding directions are the dominant  right singular vectors of the transfer function at each interpolation point. Figure \ref{fig:ln_damped_h2_evol_diff_init} shows the evolution of 
the relative $\htwo$  error during \emph{IRKA-PH} for these $5$ different selections. As before,
 in all five cases,  \emph{IRKA-PH} converges to the same reduced-model in almost the same number of steps illustrating the robustness of \emph{IRKA} to to different initializations.

\subsection{Mass-Spring-Damper system with $n=20000$}

In the previous two examples, to have a thorough analysis of the models with as many system norm computations as possible, we have chosen a very modest system size of $100$. In this example, to illustrate that we can effectively apply our method in large-scale settings, we modify the Mass-Spring-Damper model to have a full model of order $n=20000$. Then, using \emph{IRKA-PH}, we reduce the order to $r=20$, $r=30$
$r=40$ and $r=50$. In each case, the initial interpolation points are chosen as logarithmically spaced points
between $10^{-3}$ and $10^{-1}$ with the corresponding tangential directions chosen as the leading singular vectors of the transfer function at these points. Before we present the results, we note that even at this scale,
the method took less than one minute to converge with a rather straightforward implementation in Matlab. We have not tried to optimize the performance; we have simply used the Matlab sparse linear solves as is. The algorithm is expected to perform faster with the appropriate optimization of the code.

\begin{figure}[t]
\centering
\includegraphics[width=8cm]{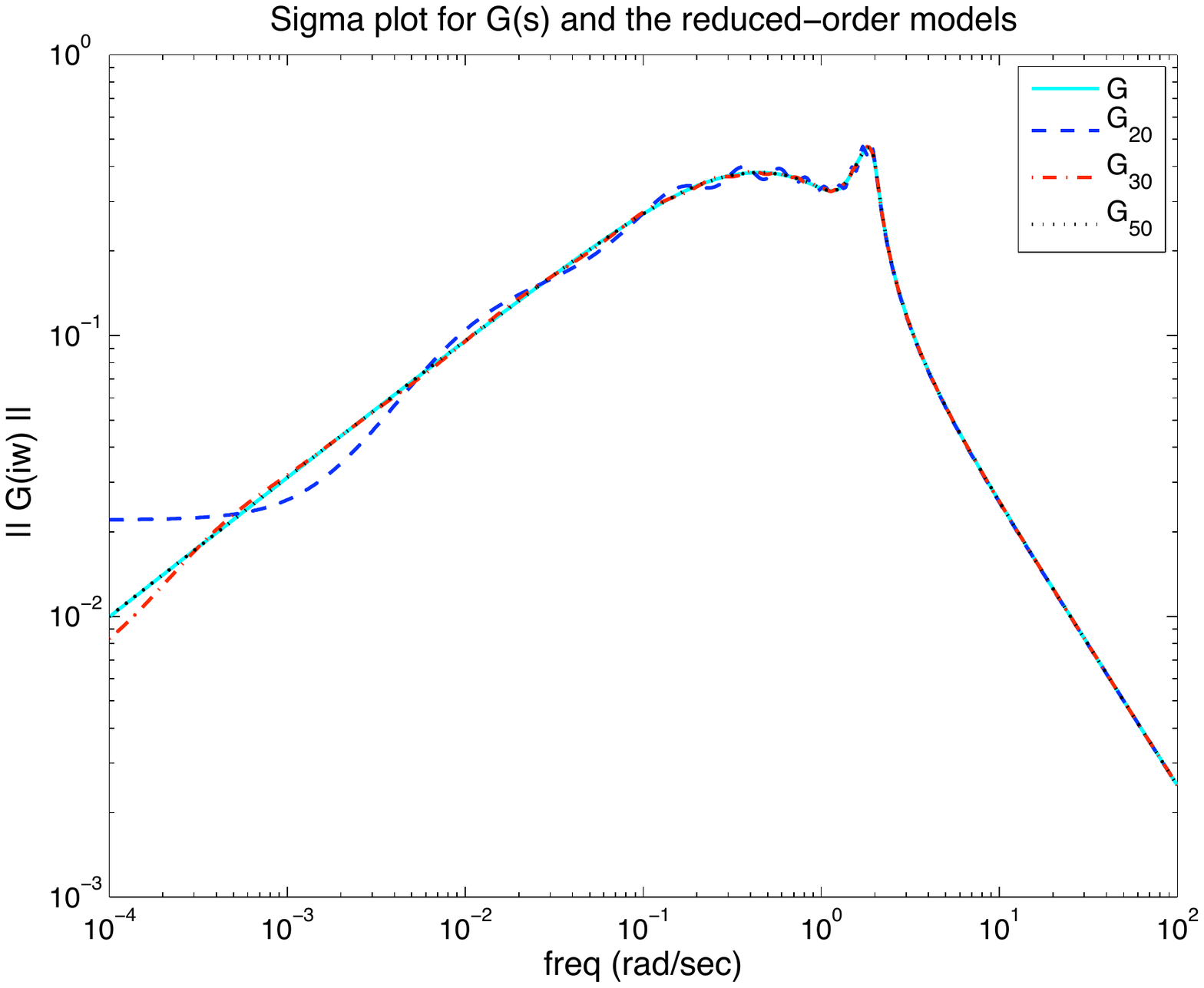}
\caption{Sigma plots of  $\bfG(s)$ and the reduced models}
\label{fig:msd_20000_sigma}
\end{figure}

\begin{figure}[t]
\centering
\includegraphics[width=8cm]{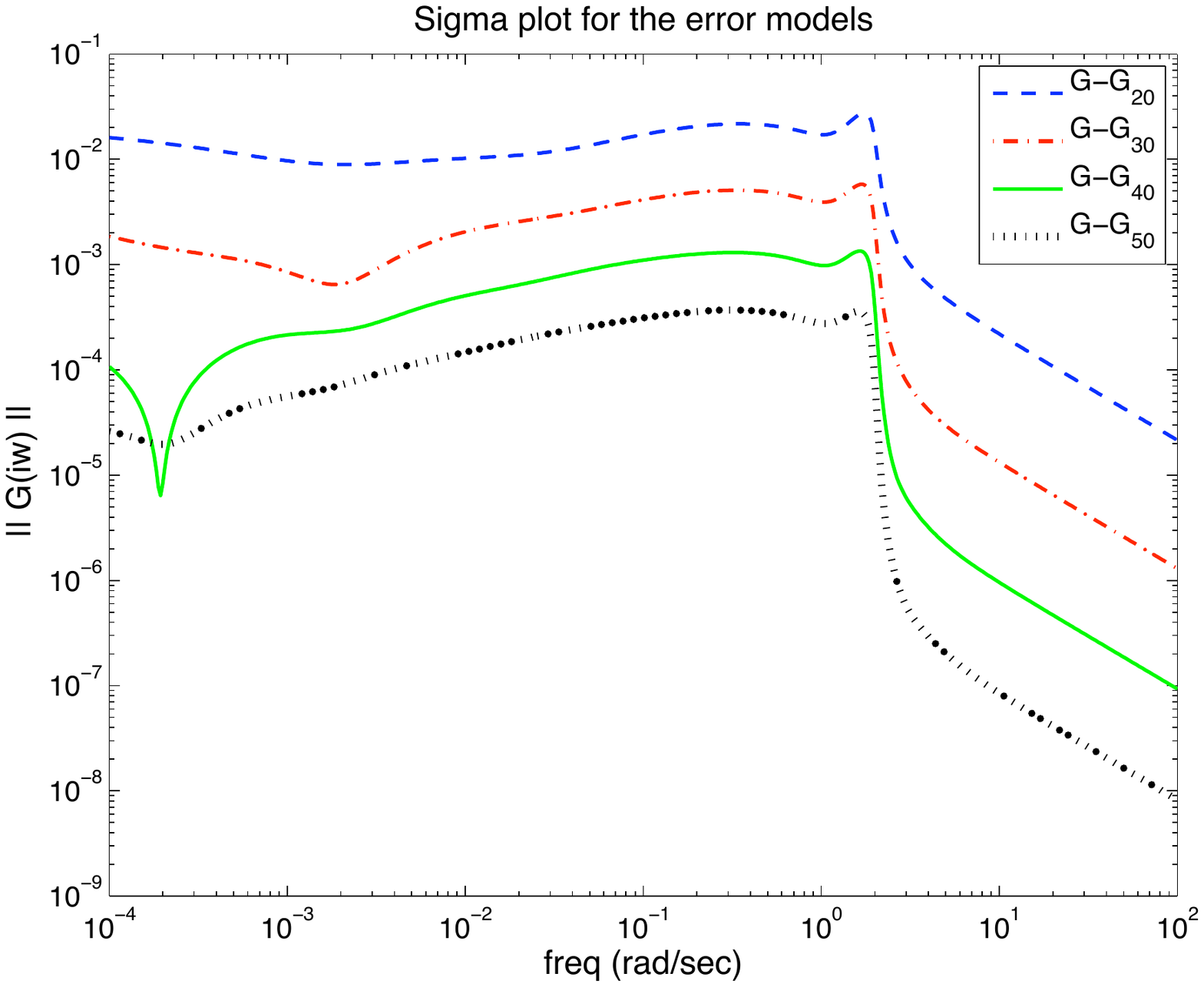}
\caption{Error sigma plots}
\label{fig:msd_20000_sigma_error}
\end{figure}

The sigma plots, i.e. $\| \bfG(\imath \omega)\|_2$ vs. $\omega$, of the full-order model $\bfG(s)$ and three of the four reduced models are plotted in Figure \ref{fig:msd_20000_sigma}.  We omitted the $40^{\rm th}$ order approximation to simplify the figure as the $r=40$ approximation was visually indistinguishable from  the $r=50$ one. Except the $r=20$ case, all the reduced models provide a high quality approximation to the full-order model of order $n=20000$. To illustrate the approximation quality further, we depict the sigma plots of the error models in Figure \ref{fig:msd_20000_sigma_error}. As $r$ increases, the quality of the approximation improves consistently; for $r=50$, we obtain an {\it approximate}\footnote{We use the term {\it approximate} since both $\| \bfG\|_\hinf$ and
$\| \bfG - \bfG_{50}\|_\hinf$ are computed by $500$ frequency sampling points over the imaginary axis as opposed to an exact $\Hinf$ norm computation. However, since the error plot is smooth, we expect this error number to be accurate enough.} relative $\Hinf$ error of $7.90\times 10^{-4}$. Hence, with the proposed algorithm, we are able to reduce a port-Hamiltonian system of order $n=20000$ in a numerically effective structure-preserving way using interpolation. Moreover, with the $\htwo$-inspired interpolation point and tangential direction selection, the reduced-model of order $r=50$ is accurate to a relative error of $7.90\times 10^{-4}$. 

\begin{figure}[t]
\centering
\includegraphics[width=8cm]{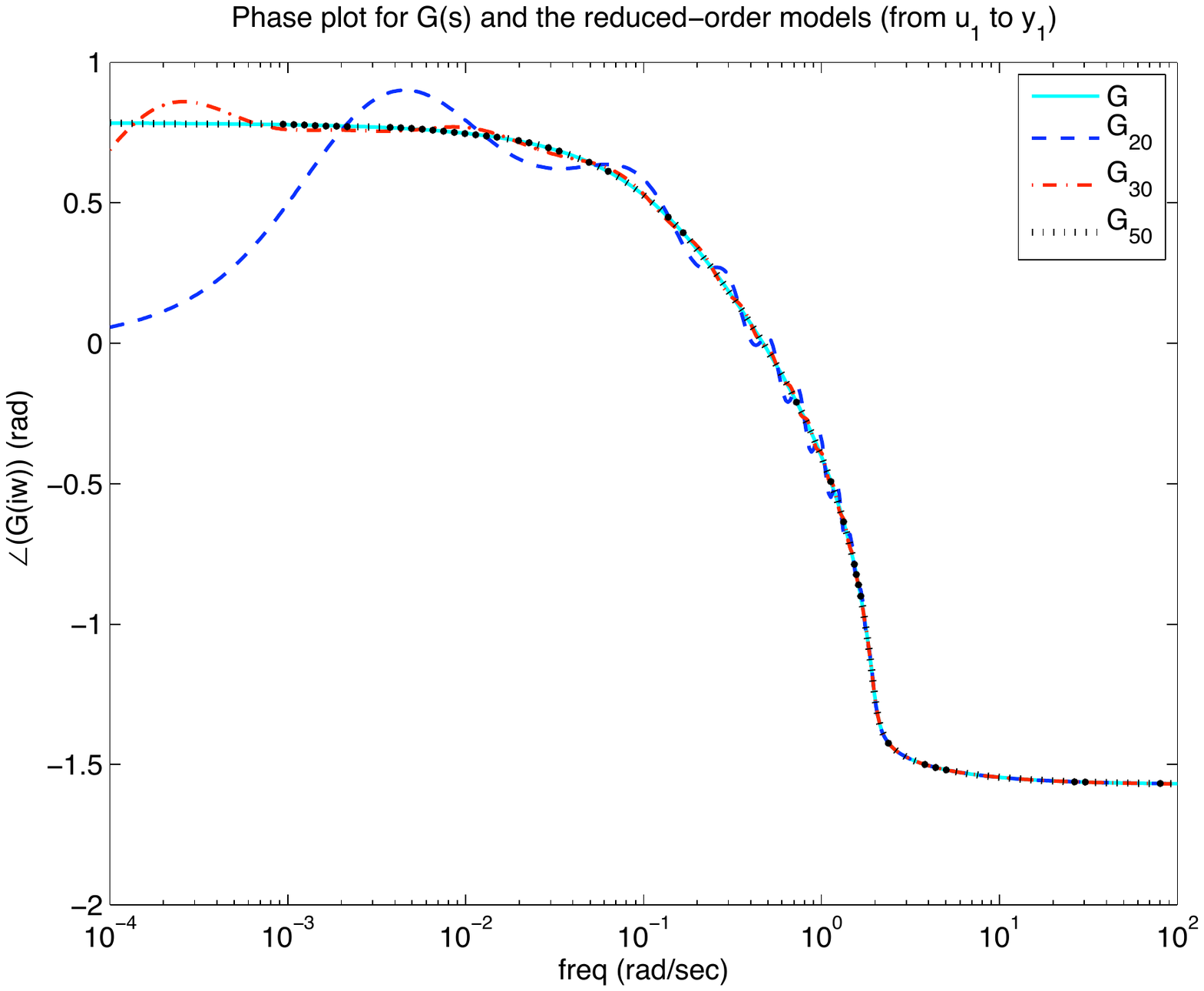}
\caption{Phase plots of  $\bfG(s)$ and the reduced models }
\label{fig:msd_20000_phase}
\end{figure}

The phase plots of the original model $\bfG(s)$ and the three of the four reduced-models are depicted in Figure \ref{fig:msd_20000_phase}.  Once again we have left the $40^{\rm th}$ order approximation out due to the same reason as before.
To make the figure simpler, we only depict the the phase plot of the subsystems relating the first input ($\bfu_1(t)$) to the first output  ($\bfy_1(t)$).  Similar to the sigma plot, Figure \ref{fig:msd_20000_phase} illustrates that the plots for the full order model $\bfG(s)$ and for the reduced model $\bfG_{40}$ and $\bfG_{50}$ are virtually indistinguishable. $\bfG_{30}$ shows some deviations around the low frequencies (as in the sigma plot). And as expected, $\bfG_{20}$ has the largest deviation. The phase plots for the other $3$ subsystems have the similar pattern and hence are omitted for the brevity of the paper.

\begin{figure}[t]
\centering
\epsfxsize = 3.3in
\epsfysize = 3.5in
\epsffile{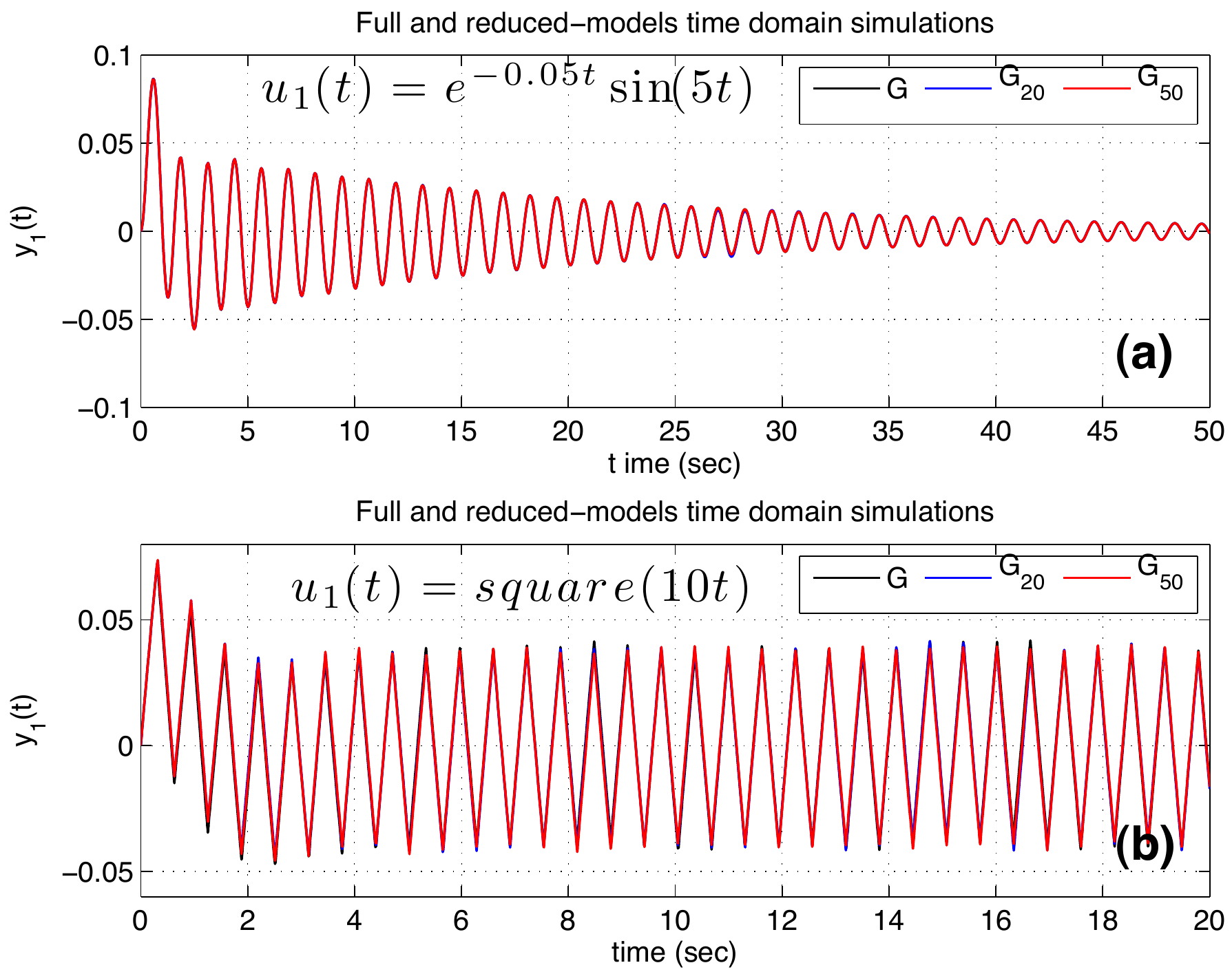}
\caption{Phase plots of  $\bfG(s)$ and the reduced models }
\label{fig:msd_20000_time}
\end{figure}

We conclude this  example by illustrating time domain simulations for the same subsystem as above 
(from the first input $\bfu_1(t)$ to the first output  $\bfy_1(t)$). We plot the results only for $\bfG(s)$, $\bfG_{20}(s)$ and $\bfG_{50}(s)$ to make the illustrations more readable.  First we use a decaying exponential $\bfu_1 =e^{-0.05 t} \sin(5t)$ and run all three models for $T=50$ seconds. The simulations are run using the \textsc{Matlab} ${\tt ode45}$ solver (taking advantage of the sparsity of the matrices in the full-order simulations).
  The results are depicted in Figure \ref{fig:msd_20000_time}-(a). Both reduced-models follow the full-order output very accurately. The maximum value of the absolute errors in the output responses are $1.32 \times 10^{-3}$ for $\bfG_{20}(s)$ and $6.15 \times 10^{-5}$ for $\bfG_{50}(s)$.  Accurate approximations are obtained in a fraction of the time, indeed. While the full-order 
simulation took $4.62$ seconds, the reduced-order simulations took $0.066$ seconds for $\bfG_{20}(s)$  and $0.068$ seconds for $\bfG_{50}(s)$; more than $98\%$ reduction in simulation time. Of course, these gains will be significantly magnified when the full-order model needs to be simulated over and over again for different input selections. The second input we try is a square wave that oscillates between $1$ and $-1$ with a period of $0.2\pi$ seconds for the time interval $0-20$ second. Figure \ref{fig:msd_20000_time}-(b) shows the corresponding outputs. Once more, both reduced models depict a high-fidelity match. For this second input,  the max value of the absolute errors in the output responses are $7.38 \times 10^{-3}$ for $\bfG_{20}(s)$ and $5.85 \times 10^{-3}$ for $\bfG_{50}(s)$. Once again, the simulation of the reduced-order models took a fraction of the time that took the full-order simulation. The simulation for $\bfG(s)$ took $9.15$ seconds. On the other hand,
it took $0.44$ second for $\bfG_{20}(s)$ and $0.46$ second for $\bfG_{50}(s)$; more than $95\%$ reduction in simulation time. As the original system dimension increases even further, these gains in simulation times would be one of the biggest advantages of model reduction.

\section{Acknowledgements}
We would like to thank the anonymous referees for their detailed reading of the manuscript and their 
valuable comments which helped improve the paper significantly.

\section{Conclusions}
We have developed a  framework for reducing multi-input/multi-output port-Hamiltonian systems via tangential rational interpolation. 
By choosing the  projection subspaces appropriately, we obtain reduced-order models
that are not only rational tangential interpolants  of the original system but they also retain the original
port-Hamiltonian structure.  Thus, they are guaranteed to be passive. 
An $\Htwo$-inspired algorithm is introduced for choosing the interpolation points and tangent directions for high-fidelity approximations.  Several numerical examples  show the effectiveness of the proposed method.
%



\bibliographystyle{plain}   
\bibliography{references}

\end{document}